\documentclass[reqno,twoside,11pt]{amsart}
\usepackage{amsmath}
\usepackage{amsfonts}
\usepackage{amssymb}
\usepackage{verbatim}
\usepackage{epsfig}

\IfFileExists{epsf.def}{\input epsf.def}{\usepackage{epsf}}

\IfFileExists{myowntimes.sty}{\usepackage{myowntimes}\usepackage{mathrsfs}}
    {\usepackage{times}\newcommand{\mathscr}{\mathcal}}
\DeclareFontFamily{OT1}{eusb}{} \DeclareFontShape{OT1}{eusb}{m}{n} {<5> <6> <7> <8> <9> <10> <11> <12> <14.4> eusb10}{}
\DeclareMathAlphabet{\eusb}{OT1}{eusb}{m}{n}

\DeclareFontFamily{OT1}{eusm}{} \DeclareFontShape{OT1}{eusm}{m}{n} {<5> <6> <7> <8> <9> <10> <11> <12> <14.4> eusm10}{}
\DeclareMathAlphabet{\eusm}{OT1}{eusm}{m}{n}

\DeclareFontFamily{OT1}{eufm}{} \DeclareFontShape{OT1}{eufm}{m}{n} {<5> <6> <7> <8> <9> <10> <11> <12> <14.4> eufm10}{}
\DeclareMathAlphabet{\mathfrak}{OT1}{eufm}{m}{n}

\DeclareFontFamily{OT1}{fraktura}{}
\DeclareFontShape{OT1}{fraktura}{m}{n} {<5> <6> <7> <8> <9> <10> <11> <12> <13> <14.4> [1.1] eufm10}{}
\DeclareMathAlphabet{\fraktura}{OT1}{fraktura}{m}{n}

\DeclareFontFamily{OT1}{cmfi}{} \DeclareFontShape{OT1}{cmfi}{m}{n} {<5> <6> <7> <8> <9> <10> <11> <12> <13> <14.4> [0.9] cmfi10}{}
\DeclareMathAlphabet{\cmfi}{OT1}{cmfi}{b}{n}

\DeclareFontFamily{OT1}{cmss}{} \DeclareFontShape{OT1}{cmss}{m}{n} {<5> <6> <7> <8> <9> <10> <11> <12> <13> <14.4> cmss10}{}
\DeclareMathAlphabet{\cmss}{OT1}{cmss}{m}{n}
\newcommand{\rrr}{{r}}

\newtheoremstyle{thm}{1.5ex}{1.5ex}{\itshape\rmfamily}{} {\bfseries\rmfamily}{}{2ex}{}

\newtheoremstyle{def}{1.5ex}{1.5ex}{\slshape\rmfamily}{} {\bfseries\rmfamily}{}{2ex}{}

\newtheoremstyle{rem}{1.3ex}{1.3ex}{\rmfamily}{} {\itshape}
{} {1.5ex}{}

\theoremstyle{thm}
\newtheorem{theorem}{Theorem}[section]
\newtheorem{lemma}[theorem]{Lemma}
\newtheorem{claim}[theorem]{Claim}
\newtheorem{proposition}[theorem]{Proposition}

\newtheorem*{Main Theorem}{Main Theorem.}
\newtheorem{corollary}[theorem]{Corollary}
\newtheorem*{special theorem}{Lindeberg-Feller Theorem for Martingales}

\theoremstyle{def}
\newtheorem{definition}[theorem]{Definition}

\theoremstyle{rem}

\numberwithin{equation}{section}

\renewcommand{\section}{\secdef\sct\sect}
\newcommand{\sct}[2][default]{%
\refstepcounter{section}
\addcontentsline{toc}{section}{{\tocsection
{}{\thesection}{\!\!\!\!#1\dotfill}}{}} \vspace{0.7cm}
\centerline{\scshape\thesection.\ #1} \nopagebreak \vspace{0.2cm}}
\newcommand{\sect}[1]{%
\vspace{0.4cm} \centerline{\large\scshape\rmfamily #1}
\vspace{0.2cm}}
\newcommand{\bepsilon}{{\bar \epsilon}}
\newcommand{\powerfour} {8}

\newcommand{\rc}{40}

\renewcommand{\subsection}{\secdef\subsct\sbsect}
\newcommand{\subsct}[2][default]{\refstepcounter{subsection}
\addcontentsline{toc}{subsection}
{{\tocsection{\!\!}{\hspace{1.2em}\thesubsection}{\!\!\!\!#1\dotfill}}{}}
\nopagebreak\vspace{0.45\baselineskip} {\flushleft\bf
\thesubsection~\bf #1.~}
\\*[3mm]\noindent
\nopagebreak}
\newcommand{\sbsect}[1]{\vspace{0.1cm}\noindent
\textbf{#1.~}\vspace{0.1cm}}

\renewcommand{\subsubsection}{%
\secdef \subsubsect\sbsbsect}
\newcommand{\subsubsect}[2][default]{%
\refstepcounter{subsubsection}
\addcontentsline{toc}{subsubsection}{{\tocsection{\!\!}
{\hspace{3.05em}\thesubsubsection}{\!\!\!\!#1\dotfill}}{}}
\nopagebreak
\vspace{0.15\baselineskip} \nopagebreak {\flushleft\rmfamily
\itshape\thesubsubsection
\ \rmfamily #1\/.}\ }
\newcommand{\sbsbsect}[1]{\vspace{0.1cm}\noindent
\rmfamily \itshape
\arabic{section}.\arabic{subsection}.\arabic{subsubsection} \
\sffamily #1\/.\ }

\renewcommand{\caption}[1]{%
\vglue0.5cm
\refstepcounter{figure}
\begin{minipage}{0.9\textwidth}\small {\sc Figure~\thefigure. }#1\end{minipage}}

\newcommand{\prob}{{\bf P}}

\newcommand{\EE}{\mathcal E}
\newcommand{\FF}{\mathcal F}

\newcommand{\KK}{\mathcal K}

\newcommand{\MM}{\mathcal M}

\newcommand{\E}{\mathbb E}

\newcommand{\N}{\mathbb N}

\newcommand{\BbbP}{\mathbb P}
\newcommand{\BbbE}{\mathbb E}

\newcommand{\R}{\mathbb R}

\newcommand{\Z}{\mathbb Z}

\newcommand{\ignore}[1]{}

\def\myffrac#1#2 in #3{\raise 2.6pt\hbox{$#3 #1$}\mkern-1.5mu\raise 0.8pt\hbox{$#3/$}\mkern-1.1mu\lower 1.5pt\hbox{$#3 #2$}}

\title[Quenched CLT for RWRE]{A quenched invariance principle for certain
ballistic random walks in i.i.d. environments}
\author{Noam Berger} \address{Department of
Mathematics,
University of California at Los Angeles,
and Department of Mathematics,
Hebrew University, Jerusalem}
\email{berger@math.huji.ac.il}
\author{Ofer Zeitouni} \address{Department of Mathematics,
University of Minnesota.}
\email{zeitouni@math.umn.edu}

\date{February 9, 2007; Revised November 20, 2007}
\begin{document}
\thanks{
%\hglue-4.5mm\fontsize{9.6}{9.6}\selectfont\copyright\,2007
%by N.~Berger and O.~Zeitouni. Reproduction, by any means, of the entire
%article for non-commercial purposes is permitted without charge.
O.Z. was partially supported by NSF grant DMS-0503775.  N.B.
 was partially supported by NSF grant DMS-0707226.} \maketitle
%Added my DMS N.

%\vspace{-5mm}
%\centerline{\textit{$^1$Department of Mathematics,
%University of California at Los Angeles,}}
%\centerline{\textit{and Department of Mathematics,
%Hebrew University, Jerusalem.}}
%\centerline{\textit{$^2$Department of Mathematics,
%University of Minnesota.}}

\begin{abstract}
We prove that every random walk in  i.i.d. environment in dimension
greater than or equal to 2  that has an almost sure positive
speed in a certain direction, an annealed invariance principle and
some mild integrability condition for regeneration times also
satisfies a quenched invariance principle. The argument is based
on intersection estimates and a theorem of Bolthausen and Sznitman.
\end{abstract}

\section{Introduction}
\label{sec:intro}\noindent

Let $d\geq 1$. A Random Walk in Random Environment (RWRE) on
$\Z^d$ is defined as follows. Let $\MM^d$ denote the space of all
probability measures on $\EE_d=\{\pm e_i\}_{i=1}^d$ and let
$\Omega=\left(\MM^d\right)^{\Z^d}$. An {\em environment} is a
point $\omega= \{\omega(x,e)\}_{{x\in\Z^d},\,{e\in\EE_d}}
\in\Omega$. Let $P$ be a probability  measure on $\Omega$. For the
purposes of this paper, we assume that $P$ is an i.i.d. measure,
i.e.
\[
P=Q^{\Z^d}
\]
for some distribution $Q$ on $\MM^d$
and that $Q$ is {\em uniformly elliptic}, i.e.
there exists a $\kappa>0$ such that for every
%$x\in\Z^d$ and
$e\in\EE_d$,
%\{\pm e_i\}_{i=1}^d$,
%\[
%P\left(\omega(x,e)<\epsilon\right)=0.
%\]
\[
Q(\{\omega(0,\cdot):\omega(0,e)<\kappa\})=0.
\]
For an environment $\omega\in\Omega$, the {\em Random Walk} on $\omega$ is a time-homogenous
Markov chain with transition kernel
\[
P_\omega\left(\left.X_{n+1}=x+e\right|X_n=x\right)=\omega(x,e).
\]
The {\bf quenched law} $P_\omega^x$ is defined to be the law on $\left(\Z^d\right)^\N$ induced by the transition kernel
$P_\omega$ and $P_\omega^x(X_0=x)=1$.
With some abuse of notation, we write $P_\omega$ also for $P_\omega^0$.
We let ${\mathcal P}^x=P\otimes P_\omega^x$ be the joint law of
the environment and the walk, and the {\bf annealed} law is
defined to be its marginal
\[
\BbbP^x(\cdot)=\int_{\Omega}P_\omega^x(\cdot)dP(\omega).
\]
We use $\E^x$ to denote expectations with respect to $\BbbP^x$.
We consistently omit the superscript $x$ if $x=0$.

We say that the RWRE $\{X(n)\}_{n\geq 0}$
satisfies the law of large numbers with
deterministic
speed $v$ if $X_n/n\to v$, $\BbbP$-a.s. For $x\geq 0$,
let
$[x]$ denote the largest integer less than or equal to $x$.
We say that the RWRE $\{X(n)\}_{n\geq 0}$
 satisfies the {\it annealed}
invariance principle with deterministic, positive definite
covariance matrix  $\sigma^2_\BbbP$
if
the linear interpolations of the  processes
\begin{equation}
    \label{eq-291206aa}
    B^n(t)=\frac{X({[nt]})-[nvt]}{\sqrt n}\,, t\geq 0
\end{equation}
converge in distribution (with respect to the supremum topology
on the space of continuous function on $[0,1]$)
as $n\to\infty$, under the measure $\BbbP$,
to a Brownian motion
of covariance $\sigma^2_\BbbP$. We say the process $\{X(n)\}_{n\geq 0}$
satisfies the
{\it quenched} invariance principle with variance $\sigma^2_\BbbP$ if for
$P$-a.e. $\omega$, the above
convergence holds under the measure $P_\omega^0$.
Our focus in this paper are conditions ensuring that when an
annealed invariance principle holds, so does a quenched one.

To state our results, we need to recall the regeneration structure
for random walk in i.i.d. environment, developed by Sznitman and
Zerner in \cite{sz-zer}. We say that $t$ is a regeneration time
(in direction $e_1$) for $\{X(\cdot)\}$  if
\[
\langle X(s),e_1 \rangle < \langle X(t),e_1 \rangle
\mbox{         whenever          }
s<t
\]
and
\[
\langle X(s),e_1 \rangle \geq \langle X(t),e_1 \rangle
\mbox{         whenever          }
s>t\,.
\]
When $\omega$ is distributed according to an i.i.d. $P$ such that
%satisfying the
%assumptions of Theorem
%\ref{thm:main}, then
the process $\{\langle X(n), e_1\rangle\}_{n\geq 0}$ is
$\BbbP$-almost surely transient to $+\infty$,
it holds by \cite{sz-zer} that, $\BbbP$-almost surely,
there exist infinitely many regeneration times for $\{X(\cdot)\}$.
Let
\[
t^{(1)}<t^{(2)}<\ldots\,,
%\quad i=1,2
\]
be all of the regeneration times for $\{X(\cdot)\}$.
Then, the sequence
$
\{(t^{(k+1)}-t^{(k)}),(X(t^{(k+1)})-X(t^{(k)}))\}_{k\geq 1}
$
is an i.i.d. sequence under $\BbbP$.
Further, if  $\lim_{n\to\infty}
n^{-1}\langle X(n),e_1\rangle>0$, $\BbbP$-a.s., then
%thanks to our assumption that the speed is positive,
we get, see \cite{sz-zer}, that
\begin{equation}\label{eq:finexp}
\E(t^{(2)}-t^{(1)})<\infty.
\end{equation}

The main result of this paper is the following:
\begin{theorem}\label{thm:main}
Let $d\geq 4$ and let $Q$ be a
uniformly elliptic distribution on $\MM^d$. Set
$P=Q^{\Z^d}$. Assume that the random walk $\{X(n)\}_{n\geq 0}$
satisfies the law of large numbers with
a positive speed in the direction $e_1$, that is
\begin{equation}
    \label{eq-291206a}
    \lim_{n\to\infty} \frac{X(n)}{n}=v\,,\BbbP-a.s\quad
    \mbox{with $v$ deterministic such that $
    \langle v,e_1\rangle>0$}\,.
\end{equation}
Assume further that the process $\{X(n)\}_{n\geq 0}$
%rescaled random walk in random environment
%$B^n(t)$, see (\ref{eq-291206aa}),
satisfies an annealed
invariance principle with variance $\sigma^2_\BbbP$.

 Assume that
%either $d\geq 5$ and
there exists an $\epsilon>0$ such that
%OOO111
    $\E(t^{(1)})^\epsilon<\infty$ and, with some $r\geq 2$,
\begin{equation}
    \label{eq-301206a}
    \E[(t^{(2)}-t^{(1)})^{r}]<\infty\,.
\end{equation}
If
$d=4$, assume further that (\ref{eq-301206a}) holds with $r>\powerfour$.
Then, the process $\{X(\cdot)\}$ satisfies a quenched invariance
principle with variance $\sigma^2_\BbbP$.
\end{theorem}
(The condition $r\geq 2$ for $d\geq 5$ can be weakened to $r>1+4/(d+4)$ by
choosing in (\ref{eq-011107h}) below $r'=r$ and modifying appropriatly
the value of $K_d$ in Proposition \ref{prop:indep} and
Corollary \ref{claim:meet_karov}.)
We suspect, in
line with Sznitman's conjecture concerning condition $T'$, see
\cite{sznitmanreview}, that (\ref{eq-301206a}) holds for $d\geq 2$ and
all $r>0$ as soon
as (\ref{eq-291206a}) holds.

A version of Theorem \ref{thm:main} for $d=2,3$ is presented in
Section \ref{addendum}. For $d=1$, the conclusion of Theorem
\ref{thm:main} does not hold, and a quenched invariance principle,
or even a CLT, requires a different centering
\cite{zeitouni,goldsheid,peterson}. (This phenomenom is typical of
dimension $d=1$, as demonstrated in \cite{RS06} in the context of the
%Another closely related one
%dimensional model where a quenched CLT requires a different
%centering is the
totally asymmetric, non-nearest neighbor, RWRE.)
Thus, some restriction on the dimension is needed.
%We suspect, in
%line with Sznitman's conjecture concerning condition $T'$, see
%\cite{sznitmanreview}, that (\ref{eq-301206a}) holds for $d\geq 2$ and
%all $r>0$ as soon
%as (\ref{eq-291206a}) holds.

%While we suspect that the critical dimension
%s $d=2$, we have only indirect evidence for this.
%
%We also suspect, in
%line with Sznitman's conjecture concerning condition $T'$, see
%\cite{sznitmanreview}, that (\ref{eq-301206a}) holds for $d\geq 2$ and
%all $r>0$ as soon
%as (\ref{eq-291206a}) holds.

Our proof of Theorem \ref{thm:main} is based on a criterion from
\cite{BoltSzni01}, which uses two independent RWRE's in the same
environment $\omega$. This approach seems limited, in principle,
to $d\geq 3$ (for technical reasons, we restrict attention
to $d\geq 4$ in the main body of the paper),
regardless of how good tail estimates
on regeneration times hold. An alternative approach to quenched
CLT's, based on martingale methods but still using the existence
of regeneration times with good tails, was developed by
Rassoul-Agha and Sepp{\"a}l{\"a}inen in \cite{RAS1}, \cite{RAS2},
and some further ongoing work of these authors. While their
approach has the potential of reducing the critical dimension to
$d=2$, at the time this paper was written,
it had not been succesful in obtaining
statements like in Theorem \ref{thm:main} without additional
structural assumptions on the RWRE.
\footnote{
%
%\noindent
%{\it
%Note added during revision:
After the first version of this paper was completed and posted,
Rassoul-Agha and Sepp{\"a}l{\"a}inen posted a preprint \cite{RS07}
in which they prove a statement similar to Theorem
\ref{thm:main}, for all dimensions $d\geq 2$, under
somewhat stronger assumptions on moments of regeneration times.
While their approach
differs significantly
from ours, and is somewhat more complicated,
we learnt from their work an extra ingredient
that allowed us to extend our approach and prove
Theorem \ref{thm:main} in all dimensions
$d\geq 2$. For the convenience of the reader, we sketch the
argument in Section \ref{addendum} below.}

Since we will consider both the case of two independent RWRE's in
different environments and the case of two RWRE's evolving in the
same environment, we introduce some notation. For $\omega_i\in
\Omega$, we let $\{X_i(n)\}_{n\geq 0}$ denote the path of the RWRE
in environment $\omega_i$, with law $P_{\omega_i}^0$. We write
$P_{\omega_1,\omega_2}$ for the law $P_{\omega_1}^0\times
P_{\omega_2}^0$ on the pair $(\{X_1(\cdot),X_2(\cdot)\})$. In
particular,
$$E_{P\times P} [P_{\omega_1,\omega_2}(\{X_1(\cdot)\}\in A_1,
\{X_2(\cdot)\}\in  A_2)]=\BbbP(\{X_1(\cdot)\}\in A_1)\cdot
\BbbP(\{X_2(\cdot)\}\in A_2)$$ represents the annealed probability
that two walks $\{X_i(\cdot)\}$, $i=1,2$, in independent
environments belong to sets $A_i$, while
%\begin{eqnarray*}
$$E_{P} [P_{\omega,\omega}(\{X_1(\cdot)\}\in A_1,
\{X_2(\cdot)\}\in  A_2)] =
%\\
 \int P_{\omega}(\{X_1(\cdot)\}\in
A_1)\cdot P_{\omega}(\{X_2(\cdot)\}\in A_2)dP(\omega)
%\end{eqnarray*}
$$
is the annealed probability for the two walks in the {\it same}
environment.

We use throughout the notation
\[
t_i^{(1)}<t_i^{(2)}<\ldots\,,\quad i=1,2
\]
for  the sequence of regeneration times of the process $\{X_i(\cdot)\}$.
Note that whenever $P$ satisfies the assumptions in Theorem
\ref{thm:main}, the estimate (\ref{eq-301206a}) holds for
$(t_i^{(2)}-t_i^{(1)})$, as well.

\noindent
{\bf Notation} Throughout, $C$ denotes a constant whose value may change
from line to line, and that may depend on $d$ and $\kappa$ only. Constants
that may depend on additional parameters will carry this dependence in
the notation. Thus, if $F$ is a fixed function  then $C_F$ denotes a constant
that may change from line to line, but that depends on $F,d$ and $\kappa$ only.
For $p\geq 1$,
$\|\cdot\|_p$ denotes the $L^p$ norm on $\R^d$ or $\Z^d$, while
$\|\cdot\|$ denotes the supremum norm on these spaces.

\section{An intersection estimate and proof of the quenched CLT}
As mentioned in the introduction, the proof of the quenched CLT involves
considering a pair of RWRE's $(X_1(\cdot),X_2(\cdot))$
in the same environment.
The main technical tool needed is the following
proposition, whose proof will be provided in Section
\ref{int_struct}.
Let  $H_K:=\{x\in\Z^d:\langle x,e_1\rangle > K\}$.
\begin{proposition}\label{prop:not_meet_insame}
    We continue under the assumptions of Theorem \ref{thm:main}.
    Let
\[
W_K:=\left\{\{X_1(\cdot)\}\cap\{X_2(\cdot)\}\cap H_K \neq\emptyset\right\}.
\]
Then
\begin{equation}\label{eq:kappa}
    E_P[P_{\omega,\omega}(W_K)]
 < C K^{-\kappa_d}
\end{equation}
where $\kappa_d=\kappa_d(\epsilon,r)>0$ for $d\geq 4$.
\end{proposition}
We can now bring the\\
{\bf Proof of Theorem \ref{thm:main}} (assuming Proposition
\ref{prop:not_meet_insame}).
For $i=1,2$, define
$B_i^n(t)=n^{-1/2}(X_i({[nt]})-[nvt]))$, where
the processes $\{X_i\}$ are RWRE's in
the same environment $\omega$, whose law is
$P$.
We
introduce
the space $C(\R_+, \R^d)$ of continuous $\R^d$-valued functions on  $\R_+$, and
the $C(\R_+, \R^d)$-valued variable
\begin{equation}\label{3.1}
\mbox{$\beta_i^n(\cdot)=$ the polygonal interpolation of $\frac{k}{n}
\rightarrow
B_i^n(\frac{k}{n})$, $k \ge 0$ }\,.
\end{equation}
It will also be useful to consider the analogously defined space $C([0,T],
\R^d)$, of continuous $\R^d$-valued functions on $[0,T]$, for $T > 0$, which we
endow with the distance
\begin{equation}\label{3.2}
d_T (v,v^\prime) = \sup\limits_{s \le T} \,|v(s) - v^\prime(s) | \wedge 1 \,.
\end{equation}
With some abuse of notation,
we continue to write $\BbbP$ for the law of the pair
$(\beta_1^n,\beta_2^n)$.
By Lemma 4.1 of \cite{BoltSzni01}, the claim will follow once we
show that for
all $T > 0$, for all bounded
Lipschitz functions $F$
on $C([0,T],\R^d)$ and $b
\in (1,2]$:
\begin{equation}\label{3.3}
\sum_m \left( E_P [E_\omega (F(\beta_1({[b^m]}))
E_\omega (F(\beta_2({[b^m]}))]-\E[
F(\beta_1({[b^m]}))]\E[
F(\beta_2({[b^m]})]
\right) < \infty \,.
\end{equation}
When proving (\ref{3.3}), we may and will assume that $F$ is bounded
by $1$ with Lipschitz constant $1$.

Fix constants $1/2>\theta>\theta'$. Write $N=[b^m]$.
Let
$$s^m_i:=\min\{t>N^{\theta}/2: \mbox{
$X_i(t)\in H_{N^{\theta'}}$,
$t$ is a regeneration time for $X_i(\cdot)$}\}\,, i=1,2.
$$
Define the events
$$A^m_i:=\{s^m_i\leq N^{\theta}\}\,,
 i=1,2,
$$
and
$$C_m:=
\{\{X_1({n+s_1^m})\}_{n\geq 0}\cap\{X_2({n+s_2^m})\}_{n\geq 0}=\emptyset\}\,,
B_m:=A^m_1\cap A^m_2\cap C_m\,.
%\{\{X^1_{t+s_1^m}\}\cap\{X^2_{t+s_2^m}\}=\emptyset\}\,.
$$
Finally, write
${\FF}_i:=\sigma(X_i(t), t\geq 0)$ and
$${\FF}_i^\Omega:=\sigma\{\omega_z: \mbox{there exists a $t$ such that
$X_i(t)=z$}\}\vee {\FF}_i, i=1,2.$$
Note that, for $i=1,2$,
\begin{eqnarray}
    \BbbP\left( (A^m_i)^c \right)&\leq &
    \BbbP(\max_{j=1}^{N^\theta} [t_i^{(j+1)}-
t_i^{(j)}]\geq N^\theta/4)+\BbbP(t_i^{(1)}>N^\theta/4)+
\BbbP(X_i(N^\theta/2)\not\in H_{N^{\theta'}})\nonumber\\
&\leq& \frac{4^{\rrr}N^{\theta} \E[(t_i^{(2)}-t_i^{(1)})^{\rrr}]}
{N^{\rrr\theta}}+
\frac{4^\epsilon\E\left([t_i^{(1)}]^\epsilon\right)}{N^{\theta\epsilon}}
+\BbbP\left(
\sum_{j=1}^{N^{\theta'}} (t_i^{(j+1)}-t_i^{(j)})>\frac{N^\theta}{4}\right)\nonumber\\
&\leq &
%Ofer-corrected noam correction
 N^{-\delta'}+4N^{\theta'-\theta}\E\left[t_i^{(2)}-t_i^{(1)}\right] \leq
2N^{-\delta'}\,,
\end{eqnarray}
%Here corrected equation to reflect lack of first moment for t_1 (N)
with $\delta'=\delta'(\epsilon,\theta)>0$ independent of $N$.
Using the last estimate and Proposition \ref{prop:not_meet_insame}, one
concludes that
\begin{equation}
    \label{eq-281206a}
    \sum_{m} E_{P} [P_{\omega,\omega}(B_m^c)]<\infty.
\end{equation}
Now,
\begin{equation}
    \label{eq-271206a}
    |\E [F(\beta_1^{[b^m]})F(\beta_2^{[b^m]})]-
    \E [{\bf 1}_{B_m}F(\beta_1^{[b^m]})
    F(\beta_2^{[b^m]})]|
\leq  \BbbP(B_m^c)\,.
\end{equation}
Let the process $\bar\beta_i^{[b^m]}(\cdot)$ be defined exactly as
the process
$\beta_i^{[b^m]}(\cdot)$, except that one replaces $X_i(\cdot)$
by $X_i(\cdot+s_i^m)$. On the event $A_i^m$, we have by construction
that
$$\sup_t \left|\beta_i^{[b^m]}(t) -
\bar \beta_i^{[b^m]}(t)\right| \leq 2N^{\theta-1/2}\,,$$
and therefore, on the event $A_1^m\cap A_2^m$,
\begin{equation}
    \label{eq-020107a}
    \left| [F(\beta_1^{[b^m]})F(\beta_2^{[b^m]})]-
     [F(\bar\beta_1^{[b^m]})
     F(\bar\beta_2^{[b^m]})]\right|\leq C N^{\theta -1/2}\,,
\end{equation}
for some constant $C$ (we used here that $F$ is Lipschitz (with constant $1$)
and bounded by $1$).

On the other hand, writing $\omega'$ for an independent copy of
$\omega$ with the same distribution $P$,
\begin{eqnarray}
    \label{eq-271206b}
    &&\E [{\bf 1}_{B_m}F(\bar\beta_1^{[b^m]})
    F(\bar\beta_2^{[b^m]})]
    =
    E_P(E_\omega
    [{\bf 1}_{B_m}F(\bar\beta_1^{[b^m]})
    F(\bar\beta_2^{[b^m]}])\nonumber\\
    &=&\E\left(
    {\bf 1}_{A_m^1}F(\bar\beta_1^{[b^m]})
    E_\omega[{\bf 1}_{A_m^2\cap C_m}
    F(\bar\beta_2^{[b^m]}])\mid {\FF}_1^\Omega]\right)\nonumber\\
    &=&E_P\left(E_\omega\left[
    {\bf 1}_{A_m^1}F(\bar\beta_1^{[b^m]})
    E_\omega[{\bf 1}_{A_m^2\cap C_m}
    F(\bar\beta_2^{[b^m]}])\mid {\FF}_1^\Omega]\right]
    \right)\nonumber\\
    &=&E_P\left(E_\omega\left[
    {\bf 1}_{A_m^1}F(\bar\beta_1^{[b^m]})
    E_{\omega'}[{\bf 1}_{A_m^2\cap C_m}
    F(\bar\beta_2^{[b^m]}])\mid {\FF}_1^\Omega]\right]
    \right)\nonumber\\
    &=&E_P\left(E_{\omega,\omega'}\left[
    {\bf 1}_{A_m^1}F(\bar\beta_1^{[b^m]})
    {\bf 1}_{A_m^2\cap C_m}
    F(\bar\beta_2^{[b^m]}])\right]\right)
    \nonumber\\
    &=&E_P\left(E_{\omega,\omega'}\left[
    {\bf 1}_{B_m}F(\bar\beta_1^{[b^m]})
    F(\bar\beta_2^{[b^m]}])\right]\right)\,.
\end{eqnarray}
%Ofer - corrected noam correction
The third equality follows from the fact that we multiply by the
indicator of the event of non-intersection.
 Since
$$
\left| E_P\left(E_{\omega,\omega'}\left[
    {\bf 1}_{B_m}F(\beta_1^{[b^m]})
    F(\beta_2^{[b^m]}])\right]\right)-
 E_P\left(E_{\omega,\omega'}\left[
    F(\beta_1^{[b^m]})
    F(\beta_2^{[b^m]}])\right]\right)\right|
    \leq \BbbP(B_m^c)\,,
    $$
    and
    $$
 E_P\left(E_{\omega,\omega'}\left[
    F(\beta_1^{[b^m]})
    F(\beta_2^{[b^m]}])\right]\right)=
 \E\left[
    F(\beta_1^{[b^m]})\right] \E\left[
    F(\beta_2^{[b^m]}])\right]
\,,$$
we conclude from the last two displays, (\ref{eq-271206b}) and
(\ref{eq-020107a}) that
%OOO111
$$  \left| \E [F(\beta_1^{[b^m]})F(\beta_2^{[b^m]})]-
 \E\left[
    F(\beta_1^{[b^m]})\right] \E\left[
    F(\beta_2^{[b^m]}])\right]\right|\leq 2\BbbP(B_m^c)+
    2CN^{\theta-1/2}\,.$$
    Together with
    (\ref{eq-281206a}), we conclude that
    (\ref{3.3}) holds, and complete the proof of Theorem
    \ref{thm:main}.
    \qed

\section{Intersection structure}\label{int_struct}
In this section we prove Proposition \ref{prop:not_meet_insame},
that is we establish estimates on the probability that two
independent walks in the same environment intersect each other in
the  half space $H_K=\{x\in\Z^d:\langle x,e_1\rangle > K\}$. It is
much easier to obtain such estimates for walks in different
environments, and the result for different environments will be
useful for the case of walks in the same environment.

% grammatical correction  N.

\ignore{ In Section \ref{sec-interd4}, we show how to extend the
analysis of this section to $d=4$, in the presence of the sronger
estimates on regeneration times assumed in Theorem \ref{thm:main}
when $d=4$.}

\subsection{The conditional random walk}
Under the assumptions of Theorem
\ref{thm:main}, the process
$\{\langle X(\cdot), e_1\rangle\}$ is $\BbbP$-a.s.
transient to $+\infty$.
Let
$$ D:=\{
\forall_{n\geq 0},\langle X(n),e_1\rangle
\geq \langle X(0),e_1\rangle\}\,.$$
By
%Together with the assumed ellipticity, we have, see
e.g. \cite{sz-zer}, we have that
\begin{equation}\label{eq:cond}
\BbbP(D)>0\,.
%\prob\left(\forall_{n\geq 0}\langle X_i(n),e_1\rangle
%\geq \langle X_i(0),e_1\rangle\right)>0
\end{equation}
%We often consider in the sequel
%the random walk $X_i$, conditioned on the event $D_i$.
%We define the conditioned random walk to be
%the random walk conditioned on the event in (\ref{eq:cond}).
%To avoid confusion, we denote by $

\subsection{Intersection of paths in
independent environments}\label{subsec:indep} In this subsection,
we let $\omega^{(1)}$ and $\omega^{(2)}$ be independent
environments, each distributed according to $P$. Let
$\{Y_1{(n)}\}$ and $\{Y_2{(n)}\}$  be random walks in the
environments (respectively) $\omega^{(1)}$ and $\omega^{(2)}$, with
starting points $U_i=Y_i(0)$. In other words, $\{Y_1{(n)}\}$ and
$\{Y_2{(n)}\}$ are  independent samples taken from the annealed
measures $\BbbP^{U_i}(\cdot)$. For $i=1,2$
%let $U_i=Y_i(0)$ be the
%starting point of the walk, and
set
$$D_i^{U_i}=
\{ \langle  Y_i(n),e_1\rangle \geq \langle  U_i,e_1\rangle \mbox{
for  } \ n\geq 0\}\,,\quad i=1,2.$$ For brevity, we drop $U_i$
from the notation and use $\prob$ for $\BbbP^{U_1}\times
\BbbP^{U_2}$ and $\prob^D$ for $\BbbP^{U_1}(\cdot|D_1^{U_1})\times
\BbbP^{U_2}(\cdot|D_2^{U_2})$.
%To simplify notation, we drop the $U_i$ from the notation and use
%$\prob$ to denote the annealed joint law of the processes $\{Y_i\}$,
%$i=1,2$.

First we prove some basic estimates. While the estimates are
similar for $d=4$ and $d\geq 5$, we will need to prove them
separately for the two cases.

\subsubsection{Basic estimates for $d\geq 5$}\label{subsubsec:5}

\begin{proposition}\label{prop:indep} ($d\geq 5$)
    With notation as above and assumptions as in Theorem
    \ref{thm:main},
\[
\prob^D\left(\{Y_1(\cdot)\}\cap\{Y_2(\cdot)\}\neq\emptyset
%\ \right| D_1^{U_1}\cap D_1^{U_2}\right)
 \right)
 < C\|U_1-U_2\|^{-K_d}
\]
%OOO111
where $K_d=\frac{d-4}{4+d}$
\end{proposition}
The proof is very similar to the proof of Lemma 5.1 of
 \cite{noam}, except that here we need a quantitative estimate
that is not needed in \cite{noam}.
\begin{proof}[Proof of Proposition \ref{prop:indep}]
We first note that the (annealed)
law of $\{Y_i(\cdot)-U_i\}$ does not depend on $i$, and is identical
to the law of $\{X(\cdot)\}$.
We also note that on the event $D_i^{U_i}$, $t_i^{(1)}=0$.
%For brevity, we use the notation $\prob$ for
%$\BbbP(\cdot|D_1^{U_1})\times \BbbP(\cdot|D_2^{U_2})$.

For $z\in\Z^d$, let
\[
F_i(z)=\prob^D(\exists_k Y_i(k)=z)
\]
and let
\[
F_i^{(R)}(z)=F(z)\cdot{\bf 1}_{\|z-U_i\|>R}.
\]
We are interested in $\|F_i\|_2$ and in $\|F_i^{(R)}\|_2$, noting
that none of the two depends on $i$ or $U_i$.
%grammatical correction. N.
We have that
\begin{equation}\label{eq:decomp}
F_i(z)=\sum_{n=1}^\infty G_i(z,n) \ \ \ \ \ \mbox{and} \ \ \ \ \
F_i^{(R)}(z)=\sum_{n=1}^\infty G_i^{(R)}(z,n)
\end{equation}
where
\[
G_i(z,n)=\prob^D(\exists_{t_i^{(n)}\leq k<t_i^{(n+1)}} Y_i(k)=z).
\]
and
\[
G_i^{(R)}(z,n)=\prob^D(\exists_{t_i^{(n)}\leq k<t_i^{(n+1)}} Y_i(k)=z)
\cdot{\bf 1}_{\|z-U_i\|>R}.
\]
are the occupation functions of $\{Y_i(\cdot)\}$.

By the triangle inequality,
\begin{equation}\label{eq:sum1}
\|F_i\|_2\leq\sum_{n=1}^\infty\|G_i(\cdot,n)\|_2
\end{equation}
and
\begin{equation}\label{eq:sum1a}
\|F_i^{(R)}\|_2\leq\sum_{n=1}^\infty\|G_i^{(R)}(\cdot,n)\|_2.
\end{equation}

Thus  we want to bound
the norm of $G_i(\cdot,n)$ and $G_i^{(R)}(\cdot,n)$.
We start with $G_i(\cdot,n)$.
Thanks to the i.i.d. structure of the regeneration slabs (see \cite{sz-zer}),
\[
G_i(\cdot,n)=Q_i^{n}\star   J\,,
\]
where $Q_i^n$ is the distribution function of $Y_i(t_i^{(n)})$ under
$\BbbP(\cdot|D_i^{U_i})$,
\[
J(z)=\prob^D(\exists_{0=t_i^{(1)}\leq k<t_i^{(2)}}Y_i(k)-U_i=z),
\]
and $\star$ denotes (discrete) convolution.
Positive speed ($\langle v,e_1\rangle>0$) tells us that
\[
\Gamma:=\|J\|_1\leq \E(t^{(2)}-t^{(1)}|D)<\infty
\]
and thus
\[
\|G_i(\cdot,n)\|_2\leq \Gamma\|Q_i^n\|_2
\]

Under the law $\prob^D$, $Q_i^n$ is the law of a sum of integrable
i.i.d. random vectors $\Delta Y_i^k=Y_i(t_i^{k+1})-Y_i(t_i^k)$,
that due to the uniform ellipticity condition are non-degenerate.
By the same computation as in \cite[Proof of claim 5.2]{noam},
we get
%, that is
%their covariance matrix is non-degenerate.
%An easy calculation (TO INCLUDE? YES!) tells us that
\[
\|Q_i^n\|_2\leq Cn^{-d/4},
\]
and thus
\begin{equation}
    \label{eq-020107c}
\|G_i^{(R)}(\cdot,n)\|_2 \leq
\|G_i(\cdot,n)\|_2\leq Cn^{-d/4}.
\end{equation}
(We note in passing that these estimates can also be obtained from a local
limit theorem applied to a truncated version of the variables $\Delta Y_i^k$.)
It follows
from the last two displays and (\ref{eq:sum1}) that for $d\geq 5$,
\begin{equation}\label{eq:fhasum}
\|F_i\|_2<C
\end{equation}

For $F_i^{(R)}$ we have a fairly primitive bound:
by Markov's inequality and the fact that the walk is a nearest neighbor walk,
for any $r'>1$,
%OOO111
%\[
%\|G_i^{(R)}(\cdot,n)\|_2\leq\|G_i^{(R)}(\cdot,n)\|_1
%\leq\BbbE({\bf 1}_{t^{(n+1)}_i>R}(t^{(n+1)}_i-t^{(n)}_i)|D)
%\leq (n+1)\frac{\E( (t^{(2)}-t^{(1)})^2|D)}{R}.
%\]
%OOO111N
\begin{eqnarray}
    \label{eq-011107h}
 &&\|G_i^{(R)}(\cdot,n)\|_2\leq\|G_i^{(R)}(\cdot,n)\|_1
 \leq\BbbE({\bf 1}_{t^{(n+1)}_i>R}(t^{(n+1)}_i-t^{(n)}_i)|D) \\
 &\leq&\BbbE\left[\right.{\bf 1}_{t^{(n)}_i>\frac R2}(t^{(n+1)}_i-t^{(n)}_i)\left|D\right]
 +\BbbE\left[\right.{\bf 1}_{t^{(n+1)}_i-t^{(n)}_i>
 \frac R2}(t^{(n+1)}_i-t^{(n)}_i)\left|D\right]\nonumber\\
 &\leq& \frac{2n\BbbE\left[t^{(n+1)}_i-t^{(n)}_i\right]}{R}
 +
C\BbbE\left(\frac{
(t^{(n+1)}_i-t^{(n)}_i)^{r'}}{R^{(r'-1)}}\right)
 %\sqrt{
 %\BbbE\left[\left.{\bf 1}_{t^{(n+1)}_i-t^{(n)}_i>\frac R2}^2\right|D\right]
 %\BbbE\left[\left.\left(t^{(n+1)}_i-t^{(n)}_i\right)^2\right|D\right]
% \nonumber\\
 \leq C\frac{n\E( (t^{(2)}-t^{(1)})^2|D)}{R}\,,
 \nonumber
\end{eqnarray}
where the choice $r'=2$ was made in deriving the last inequality.
Together with (\ref{eq-020107c}),
%and  combining with the fact that
%\[
%\|G_i^{(R)}(\cdot,n)\|_2 \leq \|G_i(\cdot,n)\|_2 \leq n^{-d/4}
%\]
we get, with $K=\left[R^{4/(d+4)}\right]$,
\begin{equation}\label{eq:arba}
%\nonumber
\|F_i^{(R)}\|_2
\leq
C\left[
\sum_{n=1}^K \frac{n}{R}
+\sum_{n=K+1}^\infty n^{-d/4}
\right]\leq
%\\
%&\leq&
C\left[
K^2/R+K^{1-d/4}
\right]
\leq
C
R^{(4-d)/(d+4)}\,.
\end{equation}

%K/R =K^{-d/4}
%K=RK^{-d/4}
%K^{1+d/4}=R
%K^{(d+4)/4}=R
%K=R^{4/(d+4)}

%R^{8/(d+4)}/R
%R^{8/(d+4)-1}
%R^{(4-d)/(d+4)}

%K^{1-d/4}
%R^{(1-d/4)4/(d+4)}
%R^{4(1-d/4)/(d+4)}
%R^{(4-d)/(d+4)}

Let $R:=\|U_2-U_1\|/2$.
An application of the Cauchy-Schwarz inequality yields
%\begin{eqnarray*}
$$\prob^D\left(\{Y_1(\cdot)\}\cap\{Y_2(\cdot)\}\neq\emptyset
%\ \right| D_1^{U_1}\cap D_2^{U_2}\right)\\
\right)
%\\
%\langle  X_i(n),e_1\rangle \geq \langle
%U_i,e_1\rangle \mbox{ for  }  i=1,2\ ;\ n\geq 0
%\right)\\
\leq
\|F_1^{(R)}\|_2^2+2\|F_1^{(R)}\|_2\|F_1\|_2
%OOO111
=O\left(R^{(4-d)/(d+4)}\right)
%\end{eqnarray*}
$$
for $d\geq 5$.
\end{proof}

Now assume that the two walks do intersect. How far from the
starting points could this happen? From (\ref{eq:arba}) we
immediately get the following corollary.
%grammatical correction. N.
\begin{corollary}\label{claim:meet_karov} ($d\geq 5$)
Fix $R$, $Y_1(\cdot)$ and $Y_2(\cdot)$ as before.
Let $A_i$ be the event that
$Y_1(\cdot)$ and $Y_2(\cdot)$ intersect, but the intersection point
closest to $U_i=Y_i(0)$ is at distance $\geq R$ from $Y_i(0)$.
Then
%OOO111
\begin{equation}\label{eq:clospoint}
    \prob^D(A_1\cap A_2)
    %|D_1^{U_1}\cap D_2^{U_2})
< C R^{(4-d)/(d+4)}.
\end{equation}
\end{corollary}

%regarding
%random walks (which are not necessarily conditional).
%{\bf EXPAND, MAKE PRECISE}

\subsubsection{Basic estimates for $d=4$}\label{subsubsec:4}

We will now see how to derive the same estimates for dimension $4$
in the presence of bounds on higher moments of the regeneration
times. The crucial observation is contained in the following lemma.
%grammatical improvement. N.
\begin{lemma}
    \label{lem-basiclemma}
    Let $d\geq 3$ and let
    $v_i$ be i.i.d., $\Z^d$-valued random variables satisfying,
    for some $r\in [2,d-1]$,
    \begin{equation}
        \label{eq-190107a}
        \langle v_1,e_1\rangle \geq 1\; \mbox{\rm a.s.} \,,\;\;\mbox{\rm and}
        \; E\|v_1\|^r<\infty\,.
    \end{equation}
    Assume that, for some $\delta>0$,
    \begin{equation}
        \label{eq-190107c}
        P(\langle v_1,e_1\rangle=1)>\delta\,,
    \end{equation}
and
    \begin{equation}
        \label{eq-190107b}
        P(v_1=z| \langle v_1,e_1\rangle=1)>\delta,\;
        \mbox{\rm for all $z\in \Z^d$ with $\|z-e_1\|_2=1$ and
        $\langle z,e_1\rangle=1$}\,.
    \end{equation}
    Then, with $W_n=\sum_{i=1}^n v_i$,
    there exists a constant $c>0$
    %and
    %$\epsilon=\epsilon(r,d)\to_{r\to\infty}\infty$
    such that for any $z\in\Z^d$,
    \begin{equation}
        \label{eq-190107d}
        P(\exists_i: W_i=z)\leq c|\langle z,e_1\rangle|^{-r(d-1)/(r+d-1)}\,,
    \end{equation}
    and, for all integer $K$,
    \begin{equation}
        \label{eq-190107e}
        \sum_{z: \langle z,e_1\rangle=K}
        P(\exists_i: W_i=z)\leq 1\,.
    \end{equation}
\end{lemma}

\begin{proof} We set $T_K=\min\{n: \langle W_n,e_1\rangle \geq K\}$.
We note first that because of
        (\ref{eq-190107c}),
for some constant $c_1=c_1(\delta)>0$ and all $t>1$,
\begin{equation}
    \label{eq-190107f}
    P(A_t)\leq
    %\#\{i\leq t: \langle v_i,e_1\rangle=1\}<c_1 t)\leq
    e^{-c_1 t}\,.
\end{equation}
    where
        $${A}_t=
        \left\{ \#\{i\leq t: \langle v_i,e_1\rangle=1\}<
        c_1 t\right\}\,.$$
Set $\bar v= E v_1$ and $v=E\langle v_1,e_1\rangle$.
%Further, by further reducing $c_1$ if needed
%       such that $c_1<E|\langle v_1,e_1\rangle|/2$,
Then, for any $\alpha\leq 1$, we get from
        (\ref{eq-190107a}) and
    the
 Marcinkiewicz-Zygmund
    inequality (see e.g. \cite[Pg. 469]{shiryayev} or, for Burkholder's
generalization,
\cite[Pg. 341]{stroock})
        that for some $c_2=c_2(r,v,\alpha)$, and all $K>0$,
        \begin{equation}
            \label{eq-190107g}
            P(T_K< K^\alpha/2v)\leq c_2 K^{-r(2-\alpha)/2}\,.
        \end{equation}
        Let ${\mathcal F}_n:=\sigma(\langle W_i,e_1\rangle, i\leq n)$
        denote the filtration generated by the $e_1$-projection
        of the random walk
        $\{W_n\}$. Denote by
        $W_n^\perp$ the projection of $W_n$
        on the hyperplane perpendicular to $e_1$. Conditioned
        on the filtration
        ${\mathcal F}_n$, $\{W_n^\perp\}$ is a
        random walk with independent (not identically
        distributed) increments, and the assumption
        (\ref{eq-190107b}) together with standard estimates
        shows that, for some constant $c_3=c_3(\delta,d)$,
        \begin{equation}
            \label{eq-190107h}
            \sup_{y\in\Z^{d-1}}{\bf 1}_{A_t^c}P(W_t=y|
            {\mathcal F}_t)\leq c_3 t^{-(d-1)/2}\,,\;
            a.s.
        \end{equation}
        Therefore, writing $z_1=\langle z,e_1\rangle$,
we get for any $\alpha\leq 1$,
%    and fixing $\alpha=2r/(r+d-1)\leq 1$,
    \begin{eqnarray}
        \label{eq-190107j}
        &&P(\exists_i: W_i=z)\leq
        P(T_{z_1}<z_1^\alpha/2v)+
        P(W_i=z\;\mbox{\rm for some $i\geq z_1^\alpha/2v$})\nonumber\\
        &\leq &
    %%OO new cor
            c_2 z_1^{-r(2-\alpha)/2 }+
            \sum_{i=z_1^\alpha/2v}^{z_1} P(W_i=z)
            \leq c_2 z_1^{-r(2-\alpha)/2 }+
            \sum_{i=z_1^\alpha/2v}^{z_1}
            E(P(W_i=z|{\mathcal F}_i))\nonumber\\
            &\leq &
            c_2 z_1^{-r(2-\alpha)/2}+
            \sum_{i=z_1^\alpha/2v}^{z_1} P(A_i)
        +   \sum_{i=z_1^\alpha/2v}^{z_1}
            E({\bf 1}_{T_{z_1}=i}
            \sup_{y\in \Z^{d-1}}
        %OOO111
            {\bf 1}_{A_i^c}P(W_i^\perp=y|{\mathcal F}_i))\nonumber\\
            &\leq &
            c_2 z_1^{-r(2-\alpha)/2}+ z_1e^{-c_1z_1^\alpha/2v}+
            c_3 (z_1/2v)^{-\alpha (d-1)/2}\sum_{i=z_1^\alpha/2v}^{z_1}
            P({T_{z_1}=i})\,,
    \end{eqnarray}
    % minor correction (N)
            where the second inequality
            uses (\ref{eq-190107g}),
    and the fifth uses (\ref{eq-190107f}) and
            (\ref{eq-190107h}).
The estimate
    (\ref{eq-190107j}) yields (\ref{eq-190107d})
by choosing
    $\alpha=2r/(r+d-1)\leq 1$.

    To see (\ref{eq-190107e}), note that the sum of probabilities
    is exactly the expected number of visits to $\{z: \langle
    z,e_1\rangle=K\}$, which is bounded by $1$.
 \end{proof}
We are now ready to state and prove the following analogue of
Proposition \ref{prop:indep}.
\begin{proposition}\label{prop:indep4}($d=4$)
    With notation as in
    Proposition \ref{prop:indep}, $d=4$ and
    $r$ in
    (\ref{eq-301206a}) satisfying
    $r>\powerfour$, we have
\[
\prob^D\left(\{Y_1(\cdot)\}\cap\{Y_2(\cdot)\}\neq\emptyset
%\ \right| D_1^{U_1}\cap D_1^{U_2}\right)
 \right)
 < C\|U_1-U_2\|^{-K_4}
\]
where $K_4>0$.
\end{proposition}
\begin{proof}
    Fix $\nu>0$ and write $U=|U_1-U_2|$. Let $\{v_i\}_{i\geq 1}$ denote an
    i.i.d. sequence of random variables, with $v_1$ distributed like
    $Y_1(t^{(2)})-Y_1(t^{(1)})$ under $\prob^D$. This sequence clearly
    satisfies
    the assumptions of
 Lemma \ref{lem-basiclemma}, with $\delta=\kappa^2 \BbbP(D)$.

 Let $T:=\BbbE(t^{(2)}-t^{(1)})$.
 By our assumption on the tails of regeneration times,
 %and
 %Chebycheff's inequality,
 %\begin{equation}
%    \label{eq-250107a}
%    \prob^D\left(t_1^{U/8T}>\frac{U}{4}\right)\leq
%    \frac{C}{U^{r/2}}\,,
% \end{equation}
% and,
for $\nu\in (0,1)$ with $\nu r>1$,
 \begin{equation}
     \label{eq-250107c}
     \prob^D\left(\exists_{i\geq \frac{U}{8T}}:
     t_1^{(i+1)}-t_1^{(i)}>i^\nu\right)\leq
     \sum_{i=U/8T}^\infty  \frac{C}{i^{\nu r}}
     \leq C U^{1-\nu r}\,.
 \end{equation}
 By Doob's maximal inequality,
 and our assumption on the tails of regeneration times,
 %and
 \begin{eqnarray}
     \label{eq-250107b}
     \prob^D\left(\exists_{i\geq \frac{U}{8T}}: t_1^{(i)}
     >2T i\right)&\leq&
     \prob^D\left(\exists_{i\geq \frac{U}{8T}}: (t_1^{(i)}-\BbbE t_1^{(i)})
     >T i\right)\nonumber \\
     &\leq&
     \sum_{j=0}^\infty
     \prob^D\left(\exists_{i\in [\frac{2^j U}{8T},
     \frac{2^{j+1} U}{8T})}: (t_1^{(i)}-\BbbE t_1^{(i)})
     >T i\right)\nonumber \\
     &\leq&
     \sum_{j=0}^\infty
     \prob^D\left(\exists_{i\leq
     \frac{2^{j+1} U}{8T}}: (t_1^{(i)}-\BbbE t_1^{(i)})
     >2^j U/8 \right)\nonumber \\
     &\leq &C \sum_{j=0}^\infty\frac{1}{(2^j U)^{r/2}}\leq
     \frac{C}{U^{r/2}}\,.
 \end{eqnarray}
 For integer $k$ and $i=1,2$, let
 $s_{k,i}=\max\{n: \langle Y_i(t_i^{(n)}),e_1\rangle\leq k\}$. Let
 $${\mathcal A}_{i,U,\nu}:=\cap_{k\geq U/8T}\{ t_i^{(s_{k,i}+1)}-
 t_i^{(s_{k,i})}\leq
 (2T k)^\nu\}\,.$$
% Changed \cup to \cap (N)
 Combining (\ref{eq-250107b}) and (\ref{eq-250107c}), we get
 \begin{equation}
     \label{eq-250107d}
     \prob^D\left( ( {\mathcal A}_{i,U,\nu})^c\right)
     \leq C [U^{1-\nu r}+U^{-r/2}]\,.\end{equation}
 For an integer $K$, set
    ${\mathcal C}_K=\{z\in \Z^d: \langle z,e_1\rangle=K\}$.
     Note that on the event
     ${\mathcal A}_{1,U,\nu}\cap
     {\mathcal A}_{2,U,\nu}$,
     if the paths $Y_1(\cdot)$ and $Y_2(\cdot)$ intersect at a point
     $z\in {\mathcal C}_K$, then there exist integers $\alpha,\beta$
     such that
     $|Y_1(t_1^{(\alpha)})-Y_2(t_2^{(\beta)})|\leq
     2(2TK)^\nu$.
    Therefore,
    with $W_n=\sum_{i=1}^n v_i$, we get from (\ref{eq-250107b}) and
    (\ref{eq-250107d}) that,
with
 $r_0=r\wedge 3$,
%For $j=1,2$, the sequence of random
% variables $\{Y_j(t_j^{i+1})-Y_j(t_j^i)\}_{i\geq 1}$ satisfy the assumptions
% in Lemma \ref{lem-basiclemma} concerning $\{v_i\}_{i\geq 1}$.
% With notation as introduced
%in the course of proving Lemma
%    \ref{lem-basiclemma}, we have
    \begin{eqnarray}
        \label{eq-190107l}
    &&\prob^D\left(\{Y_1(\cdot)\}\cap\{Y_2(\cdot)\}\neq\emptyset
    \right)\nonumber\\
    &\leq &
    %OOO111
    2\prob^D\left(t^{(U/8T)}_1\geq U/2\right)+
    %2\prob(t_1^1>\sqrt{U})+
     2\prob^D\left( ({\mathcal A}_{i,U,\nu})^c\right)
    %2\prob(\exists_i: t^{i+1}_1-t^i_1>
    %\max(U^{\nu},i^\nu))
    %+2\prob(t^{U^{1-\nu}}_1>\frac{U}{4})\nonumber\\
    %&&\quad +
    \nonumber \\
    &&+ \sum_{K>U/8T}
    \sum_{z\in {\mathcal C}_K} \sum_{z':|z-z'|<2(2TK)^\nu}
    \prob(\exists i: W_i=z)\prob(\exists j: W_j=z')\nonumber\\
    &\leq &
    C\left[U^{-\epsilon}+U^{1-\nu r}+U^{-r/2}+
%OO major
    U^{\left[1-\frac{3r_0}{r_0+3}+4\nu \right]}
    \right]\,,
\end{eqnarray}
as long as $1-3r_0/(r_0+3)+4\nu <0$,
    where Lemma \ref{lem-basiclemma}
%and the moment estimates
and
%OO minor
\eqref{eq-250107d}
%(\ref{eq-190107l})
%    on $t_2-t_1$
were used in the last inequality.
With $r>\powerfour$ (and hence $r_0=3$), one can chose $\nu>1/r$ such that all
exponents of $U$ in the last expression are negative, yielding the
conclusion.
\end{proof}

Equivalently to Corollary \ref{claim:meet_karov}, the following is
an immediate consequence of the last line of (\ref{eq-190107l})
\begin{corollary}\label{claim:meet_karov4}
    With notation as in
    Corollary \ref{claim:meet_karov},
 $d=4$ and
    $r$ in
    (\ref{eq-301206a}) satisfying
    $r>\powerfour$, we have
\[
    \prob^D(A_1\cap A_2)
    %|D_1^{U_1}\cap D_2^{U_2})
< C R^{-K^\prime_4}.
\]
with $K^\prime_4=K^\prime_4(r)>0$.

\end{corollary}

\subsubsection{Main estimate for random walks in independent
environments}\label{subsubseq:45}

Let $R>0$ and let $T^Y_i(R):=\min\{n: Y_i(n)\in H_R\}$.
%From Corollary \ref{claim:meet_karov},
%together with Comets-Zeitouni's
%epsilon coins,
%we get the following.
\begin{proposition}\label{prop:away_forever} ($d\geq 4$)
Let $Y_1(\cdot)$ and $Y_2(\cdot)$ be random walks in
independent environments satisfying the assumptions of
Theorem \ref{thm:main}, with starting points
$U_1,U_2$ satisfying $\langle U_1,e_1\rangle = \langle U_2,e_1\rangle = 0$.
Let
%OOO111
\begin{eqnarray}
    \label{eq-020107d}
    &&A(R):=\\
    &&
    \!\!\!\!\!\!\!\!\!
    \!\!\!\!\!\!\!\!\!
\{\forall_{n<T_1^Y(R)} \langle Y_1(n),e_1\rangle\geq 0\}\cap
\{\forall_{m<T_2^Y(R)} \langle Y_2(m),e_1\rangle\geq 0\}\cap
\{\forall_{n<T_1^Y(R)\,,  m<T_2^Y(R)} Y_1(n)\neq Y_2(m)\}\,.\nonumber
\end{eqnarray}
Then,
\begin{enumerate}
\item\label{item:non_meet}
There exists $\rho>0$ such that for every choice of $R$ and $U_1,U_2$ as above,
\begin{equation}\label{eq:not_meet_R}
\prob\left(A(R)\right)>\rho.
%\begin{array}{c}
%\forall_{n<T_1(R)} Y_1(n)\geq 0 \\
%\mbox{and}\\
%\forall_{m<T_2(R)} Y_2(m)\geq 0 \\
%\mbox{and}\\
%\forall_{n<T_1(R)\ ;\ m<T_2(R)} Y_1(n)\neq Y_2(m)
%\end{array}
%\right)>\rho
\end{equation}
\item\label{item:if_reg}
%Let $A(R)$ be the event in
%(\ref{eq:not_meet_R}).
Let $\hat{B}_i(n)$ be the event that
$Y_i(\cdot)$ has a regeneration time at $T_i^Y(n)$, and let
\begin{equation}
\label{eq-200907a}
 B_i(R):=\bigcup_{n=R/2}^R \hat{B}_i(n)\,.
  \end{equation}
Then
\begin{equation}\label{eq:if_reg}
\prob\left(
\left\{\{Y_1(n)\}_{n=1}^\infty\cap\{Y_2(m)\}_{m=1}^\infty\neq\emptyset\right\}
\cap A(R)\cap B_1(R)\cap B_2(R) \right)<CR^{-\beta_d}
\end{equation}
with $\beta_d=\beta_d(r,\epsilon)>0$ for $d\geq 4$.
\end{enumerate}
\end{proposition}
\begin{proof}
To see
(\ref{eq:not_meet_R}),
note first that due to uniform ellipticity,
we may and will assume that $|U_1-U_2|>C$ for a fixed arbitrary large $C$.
Since $\zeta:=\BbbP(D_1\cap D_2)>0$ does not depend on the value of $C$,
the claim then follows
%immediately from epsilon-coins
%{\bf EXPLAIN BETTER WHAT FOLLOWS}
from Propositions \ref{prop:indep} and \ref{prop:indep4} by
choosing $C$ large enough such that
$\prob^D\left(A(R)^c\right)<\zeta/2$.

To see (\ref{eq:if_reg}), note
the event $A(R)\cap B_1(R)\cap B_2(R)$ implies the event
$D_1^{U_1}\cap D_2^{U_2}$,
%, the walks
%$\{X_1(\cdot)\}$ and $\{X_2(\cdot)\}$ are conditional,
and further if
$\{Y_1(n)\}_{n=1}^\infty\cap\{Y_1(m)\}_{m=1}^\infty\neq\emptyset$
then for $i=1,2$ the closest intersection point to $U_i$ is at
distance greater than or equal to $R/2$ from $U_i$. Therefore
(\ref{eq:if_reg}) follows from Corollary \ref{claim:meet_karov}
and Corollary \ref{claim:meet_karov4}.
\end{proof}

\subsection{Intersection of paths in the same environment}\label{subsec:same}
In this subsection we take $\{X_1{(n)}\}$ and $\{X_2{(n)}\}$ to be
random walks in the same environment $\omega$, with $X_i(0)=U_i$,
$i=1,2$, and $\omega$ distributed according to $P$.  As in
subsection \ref{subsec:indep}, we also consider $\{Y_1(n)\}$ and
$\{Y_2(n)\}$, two independent random walks evolving in independent
environments, each distributed according to $P$. We continue to
use $\BbbP^{U_1,U_2}$ (or, for brevity, $\BbbP$) for the annealed
law of the pair $(X_1(\cdot),X_2(\cdot))$, and $\prob$ for the
annealed law of the pair $(Y_1(\cdot),Y_2(\cdot))$. Note that
$\prob\neq \BbbP$. Our next proposition is a standard statement,
based on coupling, that will allow us to use some of the results
from Section \ref{subsec:indep}, even when the walks evolve in the
same environment and we consider the law $\BbbP$.
%Throughout the
%entire subsection we implicitly assume the dimension-dependent
%moment estimate for the regenerations as stated in Theorem
%\ref{thm:main}.

In what follows, a stopping time $T$ with respect to the filtration
determined by a path  $X$ will be denoted $T(X)$.
\begin{proposition}\label{claim:coupling}
With notation as above,
%Let $\{X_1(\cdot)\}$ and $\{X_2(\cdot)\}$ be as above.
%Let $\{Y_1(\cdot)\}$ and $\{Y_2(\cdot)\}$ be two
%(not conditional)
%independent random walks evolving in independent environments, each
%distributed according to $P$.
%Let $\prob$
%be the annealed law for both ensembles of walks, and
let $T_i(\cdot)$, $i=1,2$
be stopping times  such that $T_i(X_i)$, $i=1,2$
are $\BbbP$-almost surely finite.
Assume $X_1(0)=Y_1(0)$
and $X_2(0)=Y_2(0)$. Set
$$I_X:=
\left\{
\{X_1(n)\}_{n=0}^{T_1(X_1)}\bigcap\{X_2(n)\}_{n=0}^{T_2(X_2)}
=\emptyset\right\} $$
 and
$$I_Y:=
 \left\{
 \{Y_1(n)\}_{n=0}^{T_1(Y_1)}\bigcap\{Y_2(n)\}_{n=0}^{T_2(Y_2)}
 =\emptyset \right\}\,.$$
Then, for any  nearest neighbor deterministic
paths $\{\lambda_i(n)\}_{n\geq 0}$, $i=1,2$,
\begin{eqnarray}\label{eq:coupling}
&&
\prob\left(
Y_i(n)=\lambda_i(n), 0\leq n\leq T_i(Y_i), i=1,2; I_Y
\right)
\nonumber
\\ &&
=
\BbbP\left(
X_i(n)=\lambda_i(n), 0\leq n\leq T_i(X_i), i=1,2; I_X
\right)\,.
\end{eqnarray}
\end{proposition}
\begin{proof}
%This is done using an easy coupling:
For every pair of
non-intersecting paths $\{\lambda_i(n)\}_{n\geq 0}$,
define three i.i.d. environments $\omega^{(1)}$, $\omega^{(2)}$
and $\omega^{(3)}$ as follows: Let
$\{J(z)\}_{
z\in\lambda_1\cup\lambda_2}$ be a collection of i.i.d. variables, of marginal
law $Q$.
 At the same time,
 let $\{\eta^{j}(z)\}_{z\in \Z^d}$, $j=1,2,3$ be three independent i.i.d.
environments, each $P$-distributed.
Then define
\[
\omega^{(1)}(z)=\left\{
\begin{array}{ll}
J(z) & \mbox{if } z\in \lambda^{(1)}\\
\eta^{(1)}(z) & \mbox{otherwise,}
\end{array}
\right.
\]
\[
\omega^{(2)}(z)=\left\{
\begin{array}{ll}
J(z) & \mbox{if } z\in \lambda^{(2)}\\
\eta^{(2)}(z) & \mbox{otherwise,}
\end{array}
\right.
\]
and
\[
\omega^{(3)}(z)=\left\{
\begin{array}{ll}
J(z) & \mbox{if } z\in \lambda^{(1)}\cup\lambda^{(2)}\\
\eta^{(3)}(z) & \mbox{otherwise,}
\end{array}
\right.
\]
and let $Y_1$ evolve in $\omega^{(1)}$, let $Y_2$ evolve in
$\omega^{(2)}$ and let $X_1$ and $X_2$ evolve in
$\omega^{(3)}$.
Then by construction,
\[
P_{\omega^{(1)},\omega^{(2)}}
 \left(
Y_i(n)=\lambda_i(n), 0\leq n\leq T_i(Y_i)
 \right)
 =
 P_{\omega^{(3)}}
 \left(
X_i(n)=\lambda_i(n), 0\leq n\leq T_i(X_i)
 \right)\,.
\]
Integrating and then summing we get (\ref{eq:coupling})\,.
\end{proof}
An immediate consequence of Proposition
\ref{claim:coupling} is that the estimates of Proposition
\ref{prop:away_forever} carry over to the processes $(X_1(\cdot),X_2(\cdot))$.
More precisely,
let $R>0$ be given and set $T^X_i(R):=\min\{n: X_i(n)\in H_R\}$.
Define $A(R)$ and $B_i(R)$ as in
    (\ref{eq-020107d}) and (\ref{eq-200907a}),
with the process $X_i$ replacing $Y_i$.
\begin{corollary}\label{cor:away_forever}
    ($d\geq 4$)
Let $X_1(\cdot)$ and $X_2(\cdot)$ be random walks in
the same environment
 satisfying the assumptions of
Theorem \ref{thm:main},
with starting points
$U_1,U_2$ satisfying $\langle U_1,e_1\rangle = \langle U_2,e_1\rangle = 0$.
%Let
%$$ A(R)=
%\{\forall_{n<T_1(R)} Y_1(n)\geq 0\}\cap
%\{\forall_{m<T_2(R)} Y_2(m)\geq 0\}\cap
%\{\forall_{n<T_1(R)\,,  m<T_2(R)} Y_1(n)\neq Y_2(m)\}\,.$$
Then,
\begin{enumerate}
\item
There exists $\rho>0$ such that for every choice of $R$ and $U_1,U_2$ as above,
\begin{equation}\label{eq:not_meet_R_X}
\BbbP\left(A(R)\right)>\rho.
%\begin{array}{c}
%\forall_{n<T_1(R)} Y_1(n)\geq 0 \\
%\mbox{and}\\
%\forall_{m<T_2(R)} Y_2(m)\geq 0 \\
%\mbox{and}\\
%\forall_{n<T_1(R)\ ;\ m<T_2(R)} Y_1(n)\neq Y_2(m)
%\end{array}
%\right)>\rho
\end{equation}
\item
%Let $A(R)$ be the event in
%(\ref{eq:not_meet_R}).
%Let $\hat{B}_i(n)$ be the event that
%$X_i(\cdot)$ has a regeneration time at $T_i^X(n)$, and let
%\[
% B_i(R):=\bigcup_{n=R/2}^R \hat{B}_i(n).
%\]
%Then,
With $C<\infty$ and $\beta_d>0$ as in (\ref{eq:if_reg}),
\begin{equation}\label{eq:if_reg_X}
\BbbP\left(
\left\{\{X_1(n)\}_{n=1}^\infty\cap\{X_2(m)\}_{m=1}^\infty\neq\emptyset\right\}
\cap A(R)\cap B_1(R)\cap B_2(R) \right)<CR^{-\beta_d}\,.
\end{equation}
\end{enumerate}
\end{corollary}
With $\beta_d$ as in (\ref{eq:if_reg_X})
and $\epsilon$ as in the statement of Theorem \ref{thm:main},
fix $0<\psi_d$ satisfying
\begin{equation}
    \label{eq-020107h}
\mbox{$\psi_d < \beta_d(1-\psi_d)$ and
$(1+\epsilon)(1-\psi_d)>1$}\,.
\end{equation}
%Let $T_i(N):=\min\{n:\langle Y_i(n),e_1\rangle \geq N\}$.
For $R$ integer, let
\begin{equation*}
 \KK_{k}(R)=\{
\exists_{(k+0.5)R^{1-\psi_d}< j < (k+1)R^{1-\psi_d}}
 \mbox{ s.t.} \, T_i(j) \, \mbox{is a regeneration time for}
 \, X_i(\cdot)
 \}\,,
\end{equation*}
 and let
\begin{equation}
C_i(R):= \bigcap_{k=1}^{\left[2R^{\psi_d}\right]}\KK_k(R).
\end{equation}
%The main goal of this section is to prove
%Proposition:
%\ref{prop:not_meet_insame}.
Proposition
\ref{prop:not_meet_insame}  will follow from the following lemma:
\begin{lemma}\label{claim:not_meet_insame} ($d\geq 4$)
    Under the assumptions of Theorem \ref{thm:main}, there
    exist  constants $C$ and $\gamma_d>0$ such that
    for all integer $K$,
\[
\BbbP\left(
%\{Y_1(\cdot)\}\cap\{Y_2(\cdot)\}\cap H_K
%\neq\emptyset \
W_K
\cap \ C_1(K)\cap C_2(K) \right)
 < C K^{-\gamma_d}\,.
\]
%where $\gamma_d>0$ for $d\geq 5$.
\end{lemma}
\begin{proof}[Proof of Lemma \ref{claim:not_meet_insame}]
Let $w:=[K^{1-\psi_d}]$ and for $k=1,\ldots,[K^{\psi_d}/2]$ define
the event
\begin{equation}\label{eq:event_H_k}
S_k=\left\{\begin{array}{cc}
\forall_{T_1(kw)\leq j<T_1((k+1)w)}X_1(j)>kw \\
{} & \mbox{and} \\
\forall_{T_2(kw)\leq j<T_2((k+1)w)}X_2(j)>kw \\
{} & \mbox{and}  \\
\{X_1(j)\}_{j=T_1(kw)}^{T_1((k+1)w)-1}\bigcap\{X_2(j)\}_{j=T_2(kw)}^{T_2((k+1)w)-1}
=\emptyset
\end{array}\right\}
\end{equation}
By
(\ref{eq:not_meet_R_X}),
\begin{equation*}
\BbbP\left(
S_k|S_1^c\cap S_2^c\cap\cdots\cap S_{k-1}^c
\right)\geq\rho\,.
\end{equation*}
Therefore,
\begin{equation}\label{eq:event_happen}
\BbbP\left(
\cup_k S_k
\right)\geq 1-(1-\rho)^{[K^{\psi_d}/2]}.
\end{equation}

Now, by
(\ref{eq:if_reg_X}),
\[
\BbbP\left(
S_k\cap C_1(K)\cap C_2(K)\cap W_K
\right) < Cw^{-\beta_d}=CK^{-\beta_d(1-\psi_d)}.
\]
We therefore get that
\[
\BbbP\left(
\cup_k S_k\cap C_1(K)\cap C_2(K)\cap W_K
\right) < CK^{-\beta_d(1-\psi_d)}K^{\psi_d}
=CK^{\psi_d-\beta_d(1-\psi_d)}.
\]
Combined with (\ref{eq:event_happen}), we get that
\[
\BbbP\left(\{X_1(\cdot)\}\cap\{X_2(\cdot)\}\cap H_K \neq\emptyset
\ \cap \ C_1(K)\cap C_2(K)
\right)
 < C K^{-\gamma_d}
\]
for every choice of $\gamma_d<\beta_d(1-\psi_d)-\psi_d$.
\end{proof}

\begin{proof}[Proof of Proposition \ref{prop:not_meet_insame}]
Note that by the moment conditions on the regeneration
times,
\[
\BbbP\left(
C_i(K)^c
\right)
\leq CK^{-\epsilon(1-\psi_d)}
+CK\cdot K^{-(1+\epsilon)(1-\psi_d)}
=CK^{-\epsilon(1-\psi_d)}+CK^{1-(1+\epsilon)(1-\psi_d)}.
\]
By the choice of $\psi_d$, see
    (\ref{eq-020107h}), it follows that
(\ref{eq:kappa}) holds for
\[
 \kappa_d<\min\left\{
 (1+\epsilon)(1-\psi_d)-1,
 \gamma_d
 \right\}.
\]
\end{proof}
%p(x>K^{1-\psi_d})=p(x^{1+\epsilon}>K^{(1+\epsilon)(1-\psi_d)})

\section{Addendum - $d=2,3$}
\label{addendum}
After the first version of this work was completed and circulated,
F. Rassoul-Agha and
T. Sepp{\"a}l{\"a}inen
have made significant
progress in their approach to the CLT, and posted an article
\cite{RS07} in which they derive the quenched CLT for all
dimensions $d\geq 2$, under a somewhat stronger assumption on the
moments of regeneration times than
(\ref{eq-301206a}). (In their work, they consider
finite range, but not necessarily nearest neighbor,
random walks, and relax the uniform ellipticity condition.) While
their approach is quite different from ours, it incorporates
a variance reduction step that, when coupled with the techniques
of this paper, allows one to extend Theorem \ref{thm:main} to all
dimensions $d\geq 2$, with a rather short proof.
In this addendum, we present the
result and sketch the proof.
\begin{theorem}\label{thm:mainadd}
    Let $d=2,3$. Let $Q$ and $\{X(n)\}$ be as in Theorem
    \ref{thm:main}, with $\epsilon=r\geq\rc$.
Then, the conclusions of Theorem \ref{thm:main} still hold.
\end{theorem}
{\bf Remark:} The main contribution to the condition $r\geq
\rc$ comes from
the fact that one needs to transfer estimates on
regenerations times in the direction
$e_1$ to regeneration times in the direction $v$, see
Lemma \ref{lem-161107} below. If $e_1=v$, or if
one is willing to assume moment bounds directly on the regeneration
times in direction $v$, then the same proof works
with $\epsilon>0$ arbitrary
and $r>14$.
\begin{proof}[Proof of Theorem \ref{thm:mainadd} (sketch)]
    The main idea of the proof is that the condition ``no
    late intersection of independent random walks in the same environment''
    may be replaced by the condition ``intersections
    of independent random walks in the same environment
    are rare''.

    Recall, c.f. the notation and proof of Theorem \ref{thm:main},
    that we need
to derive a polynomially decaying bound on $\mbox{\rm Var}(E_\omega
F(\beta^N))$ for $F:C([0,1],\R^d)\to \R$ bounded Lipschitz and
$\beta^N$ the polygonal interpolation as in (\ref{3.1}). In the
sequel, we write $F^N(X):=F(\beta^N)$ if $\beta^N$ is the polygonal
interpolation of the scaling (as in (\ref{3.1})) of the path
$\{X_n\}_{n=0,\ldots,N}$.

%added word scaling N.

For any $k$, let $S_k=\min\{n: X_n\in H_k\}$. For two paths
$p_1,p_2$ of length $T_1,T_2$ with $p_i(0)=0$,
let $p_1\circ p_2$ denote the concatenation,
i.e.
$$ p_1\circ p_2(t)=\left\{\begin{array}{ll}
    p_1(t)\,,& t\leq T_1\\
    p_1(T_1)+p_2(t-T_1)\,,& t\in (T_1,T_2]\,.
\end{array}
\right.$$
Use the notation $X_i^j=\{X_i,\ldots,X_2,\ldots,X_{j}\}$.
Then, we can write, for any $k$,
$$ F^N(X_0^N)=F^N(X_0^{S_k\wedge N}\circ[X_{S^k\wedge N }^N-
X_{S^k\wedge N}])\,.$$
Now comes the main variance reduction step, which is based on martingale
differences. Order the vertices in an $L^1$ ball of
radius $N$ centered at $0$ in $\Z^d$ in lexicographic order $\ell(\cdot)$.
Thus, $z$ is the predecessor of $z'$, denoted
$z=p(z')$, if $\ell(z')=\ell(z)+1$. Note that (because of our choice
of lexicographic order), if $z_1<z_1'$ then $\ell(z)<\ell(z')$.

Let $\delta>1/r$ be given such that $2\delta<1$. Define the event
$$ W_N:=\{\exists i\in [0,N]: t^{(i+1)}-t^{(i)}>N^{\delta}/3
\;\mbox{\rm or}\;
t^{(i+N^\delta)}-t^{(i)}>N^{3\delta/2}\}
\,.$$
By our assumptions,
 we have that
$\BbbP(W_N)\leq C(N^{-\epsilon \delta}+N^{1-\delta r})$,
and hence decays polynomially.
%Write
%OOO111
%$$ V_N:=\{\langle X_n,e_1\rangle\; \mbox{\rm does not backtrack
%more than $N^\delta/2$ steps by time $N$}\}\,.$$
%Clearly,
%$V_N^c\subset W_N$. Thus,
$$\mbox{\rm Var}(E_\omega F(\beta^N))
\leq
\mbox{\rm Var}(E_\omega F(\beta^N){\bf 1}_{W_N^c})+ O(N^{-\delta'})\,,$$
for some $\delta'>0$.
In the sequel we write $\bar F^N(X)=F^N(X){\bf 1}_{W_N^c}$.

Set
${\mathcal G}_z^N:=\sigma(\omega_x: \ell(x)\leq \ell(z), \|x\|_1\leq N)$, and
write $\hat H_k=\{z: \langle z,e_1\rangle =k\}$.
We have the following martingale difference representation:
\begin{eqnarray}
    \label{eq-cip1}
    E_\omega \bar F^N(X)-\E \bar F^N(X)&=&
\sum_{z:  |z|_1\leq N}
\left[\E \left(\bar F^N(X)|{\mathcal G}_z\right)-
\E\left(\bar F^N(X)|{\mathcal G}_{p(z)}\right)\right]\nonumber\\
&=:&\sum_{k=-N}^N
\sum_{z\in \hat H_k, |z|_1\leq N}
\Delta_z^N\,.
\end{eqnarray}
Because it is a martingale differences representation, we have
\begin{equation}
    \label{eq-cip2}
    \mbox{\rm Var}(E_\omega \bar F^N(X))=
\sum_{k=-N}^N
\sum_{z\in \hat H_k, |z|_1\leq N}
\E\left(\Delta_z^N\right)^2\,.
\end{equation}
Because of the estimate $\E[(t^{(1)})^\epsilon]<\infty$, the
Lipschitz property of $F$, and our previous remarks
concerning $W_N$, the contribution
of the terms with $k\leq 2N^\delta$  to the sum in (\ref{eq-cip2}) decays
polynomially.
To control the terms with $k>2N^{\delta}$,
for $z\in \hat H_k$ let
$\tau_z$ denote the largest regeneration time $t^{(i)}$ smaller than
$S_{{k-N^{\delta}}}$, and write $\tau_z^+$ for the
first regeneration time larger than $S_{{k+N^{\delta}}}$. Then,
$$F^N(X)=F^N(X_0^{\tau_z}\circ [X_{\tau_z}^{\tau_z^+}-X_{\tau_z}]
\circ[X_{\tau_z^+}^{N}-X_{\tau_z^+}])\,.$$
%OOO111
 Because of the Lipschitz property of $F$, our rescaling, and
 the fact that we work on the event $W_N^c$,  we
have the bound
$$|\bar F^N(X_0^{\tau_z}\circ [X_{\tau_z}^{\tau_z^+}-X_{\tau_z}]
\circ[X_{\tau_z^+}^{N}-X_{\tau_z^+}])-
\bar F^N(X_0^{\tau_z}\circ
[X_{\tau_z^+}^{N}-X_{\tau_z^+}])|\leq 4
%\max_{i=1}^N
%\frac{t^{(i+1)}-t^{(i)}}{\sqrt{N}}\leq
N^{(3\delta-1)/2}\,.
$$
One then obtains by standard manipulations
$$
\E\left(\Delta_z^N\right)^2\leq
CN^{3\delta-1}E[(E_\omega[{\bf 1}_{X\, \rm{visits }\, z}])^2]\,.$$
%Therefore, for any $p>1$ and $q=p/(p-1)$,
%$$
%\E\left(\Delta_z^N\right)^2\leq
%E[(E_\omega[{\bf 1}_{X\, \rm{visits }\, z} \frac{U_N}{\sqrt{N}}])^2]
%\leq
%[E(E_\omega[{\bf 1}_{X\, \rm{visits }\, z} ])^{2}]^{1/q}
%\frac{[\E U_N^{2p}]^{1/p}}{N}
%\,.$$
%Note that for any $r>2$ and $x\geq 1$,
%$$ \BbbP(U_N>x)\leq N \E[(t^{(2)}-t^{(1)})^r]/x^r\,.$$
%Thus,
%$$[\E U_N^{2p}]\leq 1+
% \int_{1}^\infty x^{2p-1} \BbbP(U_N>x) dx\leq
%c_rN\,$$
%for some constant $c_r<\infty$, as soon as $r>2p$.
%On the other hand, with
%$\{X_i(n)\}_{n\geq 0}$ denoting two independent
%copies of $\{X(n)\}_{n\geq 0}$ {\it in the same environment}, and
%letting
Let $I_N$ denote the number of intersections, up to time
$N$,
of  two independent
copies of $\{X(n)\}_{n\geq 0}$ {\it in the same environment}.
Then,
%we have  from H\"{o}lder's inequality that
%that for large $N$,
\begin{equation}
    \label{eq-cip3}
    \sum_{z: \|z\|_1\leq N}
    [E(E_\omega[{\bf 1}_{X\, \rm{visits }\, z} ])^{2}]
    =  E(E_{\omega\times \omega} I_N)\,.
\end{equation}
Combining these estimates, we conclude that
%for any $p>1$ with $r>2p>8$, and $q=p/(p-1)$, we have
\begin{equation}
    \label{eq-cip4}
    \mbox{\rm Var}(E_\omega \bar F^N(X))\leq CN^{3\delta-1}
%N^{(d+1)/p-1}
%   [ E(E_{\omega\times \omega} I_N)]^{1/q}
     E(E_{\omega\times \omega} I_N)
    + N^{-\delta'}\,.
%   \leq
%  N^{4/p-1}
%   [ E(E_{\omega\times \omega} I_N)]^{1/q}
%   \,.
\end{equation}
The proof of Theorem \ref{thm:mainadd} now follows from the
following lemma.
%added word lemma. N.
\begin{lemma}
    \label{lem-add}
    Under the assumptions of Theorem
    \ref{thm:mainadd}, for $d\geq 2$ and $r> \rc$, we have that for
$r'<r/4-1/2$ and any
    $\epsilon'\in (0,1/2-4/r'+2/(r')^2)$,
    \begin{equation}
        [ E(E_{\omega\times \omega} I_N)]\leq C
        N^{1-\epsilon'}\,,
    \end{equation}
    where $C$ depends only on $\epsilon'$.
\end{lemma}
Indeed,
equipped with Lemma \ref{lem-add}, we deduce from (\ref{eq-cip4}) that
$$  \mbox{\rm Var}(E_\omega F^N(X))\leq N^{-\delta'}+
CN^{1-\epsilon'}N^{3\delta-1}
\,.$$
%Thus, whenever $2d<\epsilon_{r,d}(r-2)$, we conclude that
Thus, whenever $\delta>1/r$ is chosen such that $3\delta<\epsilon'$,
(which is possible
as soon as $r>3/\epsilon'$, which in turn is possible
for some $\epsilon'<1/2-4/r'+2/(r')^2$ if $r\geq \rc$),
%$2d<\epsilon_{r,d}(r-2)$, we conclude that
$\mbox{\rm Var}(E_\omega F^N(X))\leq CN^{-\delta}$,
for some $\delta>0$. As mentioned above, this is enough to conclude.
\end{proof}
Before proving Lemma \ref{lem-add}, we need the following estimate:

\begin{lemma}\label{lem:martcalc}
Let $S_n$ be an i.i.d. random walk on $\mathbb R$ with $ES_1=0$ and
$E|S_1|^r<\infty$ for $r>3$. Let $U_n$ be a sequence of events
%such that $U_n$ is independent of $\{S_{i}-S_{i-1}\}_{i\neq n}$,
such that, for some constant $a_3>3/2$, and all $n$ large,
\begin{equation}
    \label{chic-1}
    P(U_n)\geq 1-\frac{1}{n^{a_3}}\,.
\end{equation}
In addition we assume that $\{U_k\}_{k<n}$ is independent of
 $\{S_k-S_n\}_{k\geq n}$
 for every $n$.

    Let $a_1\in (0,1)$ and $a_2>0$  be given.
Suppose further that for any $n$ finite,
   $$ P(\mbox{\rm for all $t\leq n$, $S_t\geq \lfloor t^{\frac{a_1}{2}}\rfloor$
    and $U_t$ occurs})>0\,.$$
Then, there
%We then have the following.
%    Let $a_1\in (0,1)$ and $a_2>0$  be given. Then, there
exists a constant $C=C(a_1,a_2,a_3)>0$ such that
for any $T$,
\begin{equation}
    \label{eq-1}
    P(\mbox{\rm for all $t\leq T$, $S_t\geq \lfloor t^{\frac{a_1}{2}}\rfloor$
    and $U_t$ occurs})
    \geq \frac{C}{T^{1/2+a_2}}\,.
\end{equation}
\end{lemma}

\begin{proof} Fix constants $\bepsilon>0$, $\alpha\in (0,1)$ and
$\beta\in (1,2)$ (eventually, we will take $\alpha\to 1,\beta\to 2$
and $\bepsilon\to \infty$). Throughout the proof, $C$ denote
constants that may change from line to line but may depend only on
these parameters. Define $b_i=\lfloor i^{\alpha\bepsilon}\rfloor$ and
$c_i=\lceil i^{\bepsilon+1} \rceil$. Consider the sequence of
stopping times $\tau_0=0$ and
\[
 \tau_{i+1}=\min\{n> \tau_{i}:
S_n-S_{\tau_{i}}>c_{i+1}-c_i \,\mbox{\rm or}\
S_n-S_{\tau_{i}}<b_{i+1}-c_i\}.
\]
% For $n\in [\tau_{i-1},\tau_i)$,
%set $S_n'=S_n-c_{i-1}$.
Declare an index $i$ good if
$S_{\tau_i}-S_{\tau_{i-1}}=c_{i}-c_{i-1}$.
%and $U_n$ occurs
%for all $n\in [\tau_{i-1},\tau_i)$.
Note that if the indices $i=1,\ldots,K$
are all good, then $S_n\geq b_{i-1}$ for all $n\in
(\tau_{i-1},\tau_i]$, $i=1,\ldots,K$.

Let the overshoot $O_i$ of $\{S_n\}$ at time $\tau_i$ be defined
as $S_{\tau_i}-S_{\tau_{i-1}}-(c_i-c_{i-1})$ if $i$ is good and
$S_{\tau_i}-S_{\tau_{i-1}}-(b_i-c_{i-1})$ if $i$ is not good.
By standard arguments (see e.g. \cite[Lemma 3.1]{RS07}),
$E(|O_i|^{r-1})<\infty$.
By considering the martingale $S_n$, we then get
% and the fact that
%$E(S_{\tau_i}-\cite[T1,Section 22]{spitzer} and \cite[T1,Section 23]{spitzer},
\begin{equation}
    \label{eq-2}
    P(i\, \mbox{\rm is good})\sim (1-\frac{1+\bepsilon}{i})\,,
\end{equation}
as $i\to\infty$.
By considering the martingale $S_n^2-nES_1^2$, we get
\begin{equation}
    \label{eq-3}
    E(\tau_{i+1}-\tau_i)=\Omega(i^{1+2\bepsilon})\,,
\end{equation}
as $i\to\infty$. In particular,
\begin{equation}
    \label{eq-5}
    P(\tau_{i+1}-\tau_i
    >i^{2+2\bepsilon+\delta})\leq \frac{C}{i^{1+\delta}}\,,
\end{equation}
while, from our assumption on the moments of $S_1$ and
Doob's inequality,
\begin{equation}
    \label{eq-4}
    P(\tau_{i+1}-\tau_i\leq  i^{\bepsilon \beta})
    \leq \frac{C}{i^{r\bepsilon(2-\beta)/2}}\,.
\end{equation}

\vspace{0.2cm}

We assume in the sequel that $r\bepsilon(2-\beta)/2>2$ and that
$a_3(\bepsilon \beta+1)>5+3\bepsilon+\delta$ (both these are possible
by choosing any $\beta<2$ so that $a_3\beta>3$, and then taking
$\bepsilon$ large). We say that $i+1$ is {\it very good} if it is
good and in addition, $\tau_{i+1}-\tau_i \in [i^{\bepsilon \beta},
i^{2+2\bepsilon+\delta}]$. By (\ref{eq-2}), (\ref{eq-5}) and
(\ref{eq-4}), we get
\begin{equation}
    \label{eq-6}
    P(i\, \mbox{\rm is very good})\sim (1-\frac{1+\bepsilon}{i})\,.
\end{equation}
Declare an index $i$ {\it excellent} if $i$ is very good
and in addition, $U_n$ occurs
for all $n\in [\tau_{i-1},\tau_i)$.

%We thus get that
%\begin{equation}
%    \label{eq-7}
%    P(i\, \mbox{\rm is very good for $i\leq K$})\geq \frac{C}{K^{1+\epsilon}}
%    \,.
%\end{equation}

On the event
%${\mathcal M}_K$
that
the first $K$ $i$'s are very good, we have that
%$\tau_K\leq C K^{3+2\epsilon+\delta}=:T_K$, and $S_n\geq K^{\epsilon
%\alpha}$ for $n\in [\tau_{K-1},\tau_K]$.
%On the same event,
%that the first $K$ $i$'s are very good,
%we also have
%that
$\tau_{K-1}\geq CK^{\bepsilon \beta+1}$ and $\tau_K\leq C
K^{3+2\bepsilon+\delta}=:T_K$, and $S_n\geq K^{\bepsilon \alpha}$ for
$n\in [\tau_{K-1},\tau_K]$. Letting ${\mathcal M}_K$ denote the
event that the first $K-1$ $i$'s are excellent, and $K$ is very
good, we then have, for every $n$,
\[
P(U_n^c{\bf 1}_{n\in[\tau_{K-1}, \tau_K)} |{\mathcal M}_K)\leq
K^{-a_3(\bepsilon\beta+1)}/P({\mathcal M}_K)\,.
\]

%\com{I don't understand this formula.}

%On the event ${\mathcal M}_K$
%we also have that
%$\tau_K\leq C K^{3+2\epsilon+\delta}=:T_K$, and $S_n\geq K^{\epsilon
%\alpha}$ for $n\in [\tau_{K-1},\tau_K]$.
%We thus get
%, and thus,
%conditioned on that
%event,
%$n\geq  \tau_{K}$
% implies that
%$P(U_n^c{\bf 1}_{n\geq \tau_K}
%|{\mathcal M}_K)\leq
%K^{-a_3(\epsilon\beta+1)}/P({\mathcal M}_K)\,.$

We now show inductively that
%Thus, if
$P({\mathcal M}_K)\geq C/K^{1+\bepsilon}$. Indeed,  under the above
hypotheses, we get
%$$
%P(U_n^c\, \mbox{\rm for some $n\in [\tau_{K-1},\tau_{K})$}|{\mathcal
%M}_K) \leq
%K^{-a_3(\epsilon\beta+1)+1+\epsilon+3+2\epsilon+\delta)}\,,$$
\[
 P(U_n^c\, \mbox{\rm for some $n\in [\tau_{K-1},\tau_{K})$}|{\mathcal
 M}_K) \leq
 K^{1+\bepsilon+3+2\bepsilon+\delta-a_3(\bepsilon\beta+1)}
 =K^{4+\delta-a_3+\bepsilon(3-a_3\beta)}
\]
 and
thus,
%under the above hypothesis,
with our choice of constants and
(\ref{eq-2}),
we conclude that
under the above hypothesis,
\begin{equation}
    \label{eq-6a}
    P({\mathcal M}_{K+1}|{\mathcal M}_K)
    %K\, \mbox{\rm is excellent }|
   % j\, \mbox{\rm are excellent for $j\leq K$})
    \sim (1-\frac{1+\bepsilon}{K+1})\,.
\end{equation}
%On the event that the first $K$ $i$'s are very good, we have that
%$\tau_K\leq C K^{3+2\epsilon+\delta}=:T_K$, and $S_n\geq K^{\epsilon
%\alpha}$ for $n\in [\tau_{K-1},\tau_K]$.
We thus get inductively that the hypothesis propagates and in particular
we get
\begin{equation}
    \label{eq-7}
    P(i\, \mbox{\rm is excellent for $i\leq K$})\geq \frac{C}{K^{1+\bepsilon}}
    \,.
\end{equation}
Further,
%if $T^{1/(1+\bepsilon
%\beta)}<K'<2T^{1/(1+\bepsilon \beta)}$ and the
if the first $K$ $i$'s are excellent
(an event with probability bounded below by $C/K^{1+\bepsilon}$),
%$CT^{-(1+\bepsilon)/(1+\bepsilon \beta)}$),
we have that
$\tau_{K}\leq T_K$.
Note that if $t=CK^{\bepsilon \beta+1}$ then on the above event
%(which occurs with probability bounded below by
%$Ct^{-(1+\bepsilon)/(1+\bepsilon \beta)}$),
we have that by time
$t$,
at least $Ct^{1/(3+2\bepsilon+\delta)}$ of the $\tau_i$'s are
smaller than $t$, and hence
$S_t\geq Ct^{\bepsilon \alpha/(3+2\bepsilon+\delta)}$.
We thus conclude that, for all $T$ large,
$$P(\mbox{\rm for all $t\leq T$,
$S_t\geq t^{\bepsilon \alpha/(3+2\bepsilon+\delta)}$, and $U_t$ occurs})\geq
\frac{C}{T^{(1+\bepsilon)/(1+\bepsilon \beta)}}\,.$$ Taking $\bepsilon$
large and $\beta$ close to $2$ (such that still
$r\bepsilon(2-\beta)>2$), and $\alpha$ close to $1$, completes the
proof.
\end{proof}

%\noindent
%{\bf Remark} While the lemma is stated and proved for integer-valued
%walks, the only place it was used
%was in the exit and Green function estimates
%from \cite{spitzer}.
%An easy coupling argument shows that those estimates remain valid
%without the restriction to integer valued increments.
%%remains valid
%%without this restriction.
%It is in this form that we use
%Lemma \ref{lem:martcalc} below.

\begin{proof}[Proof of Lemma \ref{lem-add} (sketch)]
Let
\[
v= \lim_{n\to\infty}\frac{X_n}{n}\neq 0
\]
%Noam Nov. 18
be the limiting direction of the random walk, and let $u$ be a unit
vector which is orthogonal to $v$.

In what follows we will switch from the regenerations in direction
$e_1$ that we used until now, and instead use regenerations in the
direction $v$, whose definition, given below, is slightly more
general than the definition of regenerations in the direction $e_1$
given in Section \ref{sec:intro}.

\begin{definition}
We say that $t$ is a regeneration time for $\{X_n\}_{n=1}^\infty$ in
direction $v$ if
\begin{itemize}
\item
$\langle X_s,v\rangle \leq \langle X_{t-1},v\rangle$ for every
$s<t-1$.
\item
$\langle X_t,v\rangle > \langle X_{t-1},v\rangle$.
\item
$\langle X_s,v\rangle \geq \langle X_{t},v\rangle$ for every $s>t$.
\end{itemize}
\end{definition}

%{\red Ofer, note that in this definition we already provide some
%minimum regeneration length, which is $\min\{\langle
%e_i,v\rangle:\langle e_i,v\rangle>0\}$ }

We denote by $t^{v,(n)}$ the succesive regeneration times of
the RWRE $X_n$ in direction $v$ (when dealing with two
RWRE's $X_i(n)$, we will use the notation $t^{v,(n)}_i$).
The sequence $t^{v,(n+1)}-t^{v,(n)}$, $n\geq 1$, is still i.i.d.,
and with $D^v$ defined in the obvious way,
the law of $t^{v,(2)}-t^{v,(1)}$ is identical to the law of $t^{v,(1)}$
conditioned on the event $D^v$.
The following lemma, of maybe independent interest,
 shows that, up to a fixed factor,
the regeneration time $t^{v,(1)}$ (and hence, also
$t^{v,(2)}-t^{v,(1)}$) inherits moment bounds
from $t^{(1)}$.
\begin{lemma}
\label{lem-161107} Assume \; $r>10$ and
 $\E((t^{(1)})^r)<\infty$. Then \
$\E(\langle X_{t^{v,(1)}}, v\rangle)^{2r'}<\infty$ \ and \
$\E((t^{v,(1)})^{r'})<\infty$ \ with $r'<r/4-1/2$.
%the regenerations $\left\{t_i^{(e_1)}\right\}$ satisfy
%$E\left[\left(t_2^{(e_1)}-t_1^{(e_1)}\right)^r\right]<\infty$ then
%$E\left[\left(t_2^{(v)}-t_1^{(v)}\right)^{r-\epsilon}\right]<\infty$.
%for every $\epsilon>0$.
\end{lemma}

\begin{proof}
On the event $(D^v)^c$, define $\tau_0=\min\{n>0: \langle X_n, v
\rangle\leq 0\}$ and
set $M=\max\{\langle X_n,v\rangle: n\in [0,\tau_0]\}$.
By \cite[Lemma 1.2]{sznitmanreview},
$\langle X_{t^{v,(1)}},v\rangle$ is
(under the annealed law)
stochastically dominated by the sum of a geometric number
of independent copies of $M+1$. Hence, if
$\E[ M^p|(D^v)^c]<\infty$ for some $p$, then
$\E|\langle X_{t^{v,(1)}}, v\rangle|^p<\infty$.

Fix a constant
 $\chi<1/2(\E t^{(1)}\wedge \E (t^{(2)}-
t^{(1)}))$ small enough so that $(2+2\|v\|_2)\chi<\|v\|_2^2$.
%Note first that by choosing $\chi$ small, on the
%event $M>x$ (all $x$ large),
%Choose $\chi$ to be small, and
Now fix some (large) number $x$. On
the event $M>x$, either
%Noam Nov. 20
\[
\begin{array}{lll}
&\bullet & t^{(\chi x)}\geq x \\
\mbox{\rm or}&&\\
 &\bullet & t^{(k+1)}-t^{(k)}\geq
\chi k\,\quad\mbox{\rm for \  some $k>\chi x$}\\
\mbox{\rm
or}&&\\
&\bullet &
%\begin{array}{l}
%\{t^{(k+1)}-t^{(k)}<\chi k\}\,, \mbox{\rm all}
%\; k\geq \chi x\\
%\mbox{\rm and}\\
\{|t^{(k)}-\E t^{(k)}|>\chi k\} \;\mbox{\rm or} \;
\{\|X_{t^{(k)}}-\E X_{t^{(k)}}\|>\chi k\}\,\; \mbox{\rm for \
some}\, k>\chi x\,.
%\end{array}
\end{array}
\]
(Indeed, on the event $M>x$ with $x$ large, the RWRE
 has to satisfy that
 at some large time $t>x$,
 $\langle X_t,v \rangle$ is close to $0$ instead of close to $\|v\|_2^2
 t$.)
 % More precisely, $|X_t-vt|\geq\langle
 %X_t-vt,v\rangle>\langle vt,v\rangle=\|v\|^2_2t$.
 %Let $k$ be such that $t^{(k)}$ is the last $e_1$ regeneration
 %before that time $t$. Then either $k<\chi x$, or $t-t^{(k)}>\chi k$
 %k$)
%Noam Nov. 18

%{\red We might want to multiply the last $\chi k$ by $|v|$}

Due to the moment bounds on $t^{(1)}$ and $t^{(2)}-t^{(1)}$,
and the chosen value of $\chi$,
we have
$\BbbP ( t^{(\chi x)}\geq x)\leq C x^{-r/2}$.
We also have
%$$
% \BbbP (  t^{(k+1)}-t^{(k)}\geq
%\chi k\,,\quad\mbox{\rm some $k>\chi x$})\leq C
%\sum_{k=\chi x}^\infty k^{-r/2} =C x^{-r/2+1}\,,$$
$$
 \BbbP (  t^{(k+1)}-t^{(k)}\geq
\chi k\,,\quad\mbox{\rm some $k>\chi x$})\leq C \sum_{k=\chi
x}^\infty k^{-r} =C x^{-r+1}\,,$$ and
%Noam Nov. 18
$$\BbbP( \{|t^{(k)}-\E t^{(k)}|>\chi k\}
\;\mbox{\rm or}
\;
\{\|X_{t^{(k)}}-\E X_{t^{(k)}}\|>\chi k\}\,,\;
\mbox{\rm some}\, k>\chi x)
\leq C \sum_{k=\chi x}^\infty k^{-r/2} =C x^{-r/2+1}\,.$$
%\end{array}
We conclude that $\E M^p\leq C+C\int_1^\infty x^{p-1} x^{-r/2+1}dx
<\infty$ if $p<r/2-1$. This proves that
\begin{equation}
\label{eq-171107a}
\E|\langle X_{t^{v,(1)}}, v\rangle|^p<\infty\quad \mbox{\rm if $p<r/2-1$}.
\end{equation}

We can now derive
moment bounds on $t^{v,(1)}$ (which imply also
moment bounds on $\bar t^v:=t^{v,(2)}-t^{v,(1)}$).
Clearly, $E(\bar t^v)<\infty$, and it is easy to see that this fact and
the higher moments 
bounds on regeneration times in the $e_1$ direction (third moment assumption
suffices) imply that also
$E(t^{v,(1)})<\infty$.
%Although it is not written, it is true (without conditioning) by running
%exactly the same argument as in the given proof. Explicitely: let
%me write $t^v$ for $t^{(v),1}$. Suppose
%$E(t^v)=\infty$. Now,
%this implies
%$$P(A_1):=P(t^{v}>x_m)\geq x_m^{-1-\epsilon}$$ for an appropriate
%sequence of $x_m$. But you already know (by the assumed moment bounds on
%$t^1$) that for $t$ large,
%$P(|X_t-tv|>\delta t)\leq 1/t^3$, say.
%Thus,
%$$P(A_2):=P(\mbox{\rm for some $t>x_m$}, |X_t-tv|>\delta t)\leq 1/x_m^2$$
%In particular,
%$$P(A_1\cap A_2^c)\geq 1/2x_m^{1+\epsilon}$$
%This tells you that
%$P(\langle X_{t^{v}}, v\rangle >\delta' x_m)\geq 1/2x_m^{1+\epsilon}\,.$
%But this contradicts the moment bounds on the regeneration distance.
%
Suppose
$\E( (t^{v,(1)})^{p'})=\infty$.
% but
%$\E (\bar t^v)^{q}<\infty$ for $q<p'$.
For any $\epsilon''>0$
we can then
find a sequence of integers
$x_m\to\infty$
such that
$\BbbP(t^{v,(1)}>x_m)\geq C/x_m^{p'+\epsilon''}$.
Therefore, using (\ref{eq-171107a}) and the assumed moment bounds,
\begin{eqnarray*}
&& \BbbP(|t^{v,(x_m)}-\E (t^{v,(x_m)})|> x_m/2)
\\
&\geq& \BbbP(t^{v,(1)}-\E(t^{v,(1)})>x_m)\BbbP(|t^{v,(x_m)}-t^{v,(1)}-
(x_m-1)\E (\bar t^v)|<\chi x_m)
\geq Cx_m^{-(p'+\epsilon'')}\nonumber \,.
\end{eqnarray*}
%Noam Nov. 18
Therefore,
%Ofer Nov 20 modified epsilon''
\begin{eqnarray}\label{eq-171107b}
 \nonumber
 \BbbP\left(
 \left|t^{v,(x_m)}-
 \E\left[t^{v,(x_m)}\right]\right|>x_m/2
 \ \ ; \ \
 \left|\langle X_{t^{v,x_m}}, v\rangle-
 \E \langle X_{t^{v,x_m}}, v\rangle\right|<\chi x_m
 \right)\\
\geq
\frac{C}{x_m^{p'+\epsilon''}}-\frac{C}{x_m^{p/2}} \geq
\frac{C}{x_m^{p'+\epsilon''}}\,,
\end{eqnarray}
%Noam Nov. 18
if $p'<p/2<r/4-1/2$. On the other hand, the event depicted in
(\ref{eq-171107b}) implies that at some time $t$ larger than $x_m$,
the ratio $\langle X_t, v\rangle/t$ is not close to $\|v\|_2^2$, an
event whose probability is bounded above (using the regeneration
times $t^{(n)}$) by
$$Cx_m^{-r/2}+C\sum_{k=Cx_m}^\infty k^{-r/2}\leq Cx_m^{1-r/2}\,.$$
Since $1-r/2<-p'$, we achieved a contradiction.
%
%
%First, we show that
%\begin{equation}\label{eq:dist}
%E\left[\left\langle
%X_{t_2^{(v)}}-X_{t_1^{(v)}},v\right\rangle^{r-\epsilon}\right]<\infty
%\end{equation}
%\newcommand{\hx}{{\bar{X}}}
%For ease of notation, we denote $\bar{x}:=\langle x,v\rangle$.
% Without loss of generality, we condition on
%the event $t_1=0$.
%
% As in \cite{sz-zer}, we let $h_1,h_2,\ldots$ be
%the successive maxima in the $v$ direction, and we remember that
%$\hx_{t_2^{(v)}}\leq h_T+1$ with $T$ being geometric and independent
%of $\{h_j:j\leq T\}$. Note also that if $j\leq T$ then after we
%reach distance $h_j$, we retreat to distance $h_{j-1}$. Therefore in
%order to prove (\ref{eq:dist}) it is sufficient to show that
%\[
%\sup\{E\left[(h_{j+1}-h_j)^r
%%{\bf 1}_{j\leq T}
%\right]:j=1,2,\ldots\}<\infty
%\]
%
% For given $K$ and $j$, we want to estimate the
%probability $\prob\left(h_{j+1}-h_j>K\right)$.
%
%Let $\tau_j$ be the time at which $h_j$ is arrived at.
%
%Let $\nu=X_{\tau_{j+1}}-X_{\tau_{j+1}}$.
%
%Let $b=\lim_{n\to\infty}\frac{\langle X_n,e_1\rangle}{n}$.
%Then
%\[
% \prob(h_{j+1}-h_j>K) \leq
% =\prob\left(\langle \nu,j\rangle>K \mbox{ and } \langle \nu,e_1\rangle\leq\frac{bK}2\right)
% +\prob\left(\langle \nu,e_1\rangle>\frac{bK}2\right)
%\]
%
%Now,
%\[
%\prob\left(\langle \nu,e_1\rangle>\frac{bK}2\right)<K^{-r}
%\]
%because of the moment estimates for the $e_1$ regenerations.
%
%To bound
%\[
%\prob\left(\langle \nu,j\rangle>K \mbox{ and } \langle
%\nu,e_1\rangle\leq\frac{bK}2\right)
%\]
%
%we note that {\bf here there is some more detailed structure to deal
%with}
%
\end{proof}

Consider temporarily the walks $X_1$ and $X_2$ as evolving in
independent environments. We define the following i.i.d. one
dimensional  random walk:
\[
\begin{array}{c}
%W_n=\left\langle X_1\left(t_1^{(n)}\right),v \right\rangle \\
S_n=\left\langle
X_1\left(t_1^{v,(n)}\right)-X_2\left(t_2^{v,(n)}\right),u \right\rangle.
\end{array}
\]
\newcommand{\tempor}{\eta}
\newcommand{\horiz}{\kappa}
\newcommand{\horiza}{\kappa}
\newcommand{\horizb}{\eta}
Set $r'<r/4-1/2$.
For $\horiza$ and $\horizb$ to be determined below, we define the
events
 $$B_n=\left\{
 t_i^{v,(n)}-t_{i}^{v,(n-1)}<n^{\tempor}, i=1,2\right\}\,,$$
\[
 C_n=
 \left\{
\left| \langle
 X_i(t_i^{v,(n)})-E(X_i(t_i^{v,(n)})),v
 \rangle \right|
 <n^{\horiz}, i=1,2
 \right\}\,,
\]
\[
 D_n=
 \left\{
 \max_{k\in[n-n^{\horiza},n]}
 %\stackrel{k\leq n}{|k-n|\leq n^{\horiza}}}
 |\langle
 X_i(t_i^{v,(k)})-X_i(t_i^{v,(n)}),u
 \rangle|
 <n^{\horizb}, i=1,2
 \right\}\,,
\]
%\com{\bf Are $\alpha$ and $\beta$ as in Lemma \ref{lem:martcalc}?}
and  $U_n=B_n\cap C_n\cap D_n$. By our assumptions,
Lemma \ref{lem-161107},
 and standard
random walk estimates, $P(B_n^c)\leq n^{-\tempor r'}$, $P(C_n^c)\leq
n^{-r'(2\horiza-1)}$ and $P(D_n^c)\leq n^{-r'(2\horizb-\horiza)/2}$.
With $r'>15/2$, choose $\horiza>1/2, \horizb<1/2$ such that $r'\tempor>3/2$,
$r'(2\horiza-1)>3/2$ and $r'(2\horizb-\horiza)>3$, to deduce that
$P(U_n^c)\leq n^{-a_3}$ for some $a_3>3/2$. (This is possible with $\eta$
close to $1/2$ and $\kappa$ close to $2(\eta+1)/5$.)

%Fix $\upsilon>0$ and define the event
Fix $\eta'\in(0,1/2-\eta)$ and
define the event
$$A(T)=\{
    \mbox{\rm for all $n\leq T$,
    $S_n\geq \lfloor n^{\frac{1}{2}-\eta'}\rfloor\}$}\,.$$
    Note that there exists $k_0$ such that
    on the event $A(T)\cap_{n=1}^T U_n$,
$X_1[k_0,T/2]\cap X_2[k_0,T/2]=\emptyset$.

%While $\{S_n\}$ is not integer valued,
% it does take values
%on a multiple of $\Z$, and thus
%as mentioned above
From Lemma
\ref{lem:martcalc} we have that
%can still be applied  to yield that
    $\BbbP^D(A(T)\cap_{n=1}^T U_n)
    \geq {C}/{T^{1/2+a_2}}\,,$
    for some constant $a_2>0$. By ellipticity, this implies
\begin{equation}
\label{eq-231007a}
\BbbP^D(
    X_1[1,t_1^{(T)}]\cap X_2[1,t_2^{(T)}]=\emptyset)
    \geq {C}/{T^{1/2+a_2}}\,.
\end{equation}
    uniformly over the starting points.
%    From this, it is easy to deduce the required estimate
%    on the number of intersection. We omit further details.
(This estimate, which was derived initially for walks in independent
environments, obviously holds for walks in the same environment, i.e.
under $\BbbP$, too, because it involves a non-intersection event.)

Fix $T$, and let
\[
G(T)=\sum_{i,j}{\bf 1}_{X_1(i)=X_2(j)}{\bf 1}_{\langle
X_1(i),v\rangle\in[T-0.5,T+0.5)}.
\]
We want to bound the sum of
\[
F(T)=E\left[E_{\omega\times\omega}\left(G(T)\right)\right].
\]

%First, we get an estimate for
%$E\left[P_{\omega\times\omega}\left(G(T)\neq 0\right)\right]$:

%A renewal argument combined with (\ref{eq-231007a}) and
%an estimate on the size of the
%largest among the first $T$ regenerations tells us that
%\begin{equation}
%\label{eq:ren2}
%\sum_{t=1}^N E\left[P_{\omega\times\omega}\left(G(t)\neq
%0\right)\right]\leq CN^{1/2+a_2+1/r}
%\,.
%\end{equation}
%and therefore
%\begin{equation}\label{eq:ren2}
%\sum_{t=1}^N \sqrt{E\left[P_{\omega\times\omega}\left(G(t)\neq
%0\right)\right]}\leq CN^{\frac 34 + \frac {a_2}{2} + \frac {1}{2r}}
%\,.\end{equation}
\begin{claim}\label{lem:ren}
\begin{equation}
\label{eq:ren2} \sum_{t=1}^N
E\left[P_{\omega\times\omega}\left(G(t)\neq 0\right)\right]\leq
CN^{1/2+a_2+1/r'} \,.
\end{equation}
%and therefore
%\begin{equation}\label{eq:ren2}
%\sum_{t=1}^N \sqrt{E\left[P_{\omega\times\omega}\left(G(t)\neq
%0\right)\right]}\leq CN^{\frac 34 + \frac {a_2}{2} + \frac {1}{2r}}
%\,.\end{equation}
\end{claim}

\begin{proof}
We define variables $\{\psi_n\}_{n=1}^\infty$ and
$\{\theta_n\}_{n=0}^\infty$ inductively as follows:
\[\psi_1:=\max\{t_1^{v,(1)},t_2^{v,(1)}\}\,,\quad \theta_0=0\]
and then, for every $n\geq 1$,
% Ofer November 20 changed tau^{(i)} to $\tau_i$.
\[
\theta_n:=\min\{k>\psi_n:G(k)\neq 0\}\,,\quad
 \psi_{n+1}:=\max\{
 \tau_{1}(\theta_n),\tau_{2}(\theta_n)
 \}\,,
\]
%and $\psi_{n+1}$ is defined to be
%\[
% \psi_{n+1}:=\max\{
% \tau^{(1)}(\theta_n),\tau^{(2)}(\theta_n)
% \}
%\]
with
\[
 \tau_{i}(k) :=
 \min\{
 \langle X_{t_i^{v,(m)}},v\rangle
 :
 \langle X_{t_i^{v,(m)}},v\rangle>k+1
 \}.
\]

We define $h_n=\psi_n-\theta_{n-1}$  and
$j_n=\theta_n-\psi_n$. By (\ref{eq-231007a}), for every $k$
\begin{equation}\label{eq:momjn}
\BbbP(j_n>k|j_1,\ldots,j_{n-1},h_1,\ldots,h_n) \geq
{C}/{k^{1/2+a_2}}\,.
\end{equation}
Let
\[
K:=\min\left\{n:\sum_{i=1}^n\ j_i>N\right\}.
\]
Let
$Y_i^{(N)}=\max_{k=0}^N [t^{v,(k+1)}_i-t_i^{v,(k)}]$  be the
length of the longest of the first $N$ regenerations of $X_i$,
$i=1,2$, in direction $v$, and set $Y_N=\max(Y_1^{(N)},Y_2^{(N)})$.
Then
\[
\sum_{t=1}^N{\bf 1}_{G(t)\neq 0}\leq K\cdot Y_N
\,.\]
% with
%\[
%Y_n:=\max\left\{t_n^{(i)}-t_{n-1}^{(i)} : n\leq N\ ;\ i=1,2\right\}
%\]
We see below in (\ref{eq:yomtov}) that $E(Y_N^p)\leq CN^{p/r'}$ for
$p<r'$. In addition, by the moment bound (\ref{eq:momjn}), for any
$t$,
\[
\prob(K>t) =\prob\left(\sum_{i=1}^t j_i < N\right)
 \leq\exp\left(-C\frac{t}{N^{\frac 12+a_2}}\right).
\]
From here we get
\begin{equation}\label{eq:epstagtag}
\sum_{t=1}^N E\left[P_{\omega\times\omega}\left(G(t)\neq
0\right)\right]\leq CN^{1/2+a_2+1/r'+\epsilon''} \,
\end{equation}
for every $\epsilon''>0$. The fact that $a_2$ was an arbitrary
positive number allows the removal of $\epsilon''$ from
(\ref{eq:epstagtag}).
\end{proof}

In addition, $G(t)$ is bounded by the product of the length of the
$\{X_1\}$ regeneration containing $t$ and that of the $\{X_2\}$
regeneration containing $t$. So
%if $Y_i^{(N)}$ is defined to be the
%length of the longest of the first $N$ regenerations of $X_i$,
%$i=1,2$, then
for all $t<N$,
\[
G(t)\leq Y_1^{(N)}\cdot Y_2^{(N)},
\]
and therefore, for any $p<r'/2$,
\[
E\left[G(t)^p\right]\leq
\sqrt{E((Y_1^{(N)})^{2p})E((Y_2^{(N)})^{2p})}\leq CN^{2p/r'} \,,\]
where in the last inequality we used the estimate
\begin{equation}\label{eq:yomtov}
E( (Y_i^{(N)})^{2p}) \leq A^{2p}+ 2pN\int_A^\infty
y^{2p-1}P(\tau_i^{(2)}-\tau_i^{(1)}>y) dy\leq A^{2p}+CNA^{2p-r'}\,,
\end{equation}
with $A=N^{1/r'}$. Thus, with $1/q=(p-1)/p$,
%\begin{eqnarray*}
\begin{eqnarray}
\label{eq-231007b}
 E[G(t)]&=&E[G(t)\cdot {\bf 1}_{G(t)\neq 0}]
 \leq\left(EG(t)^p\right)^{1/p}
 \left(E\left[P_{\omega\times\omega}\left(G(t)\neq 0\right)\right]\right)^{1/q}
\nonumber\\
&\leq &
 CN^{2/r'}
\left(E\left[P_{\omega\times\omega}\left(G(t)\neq 0\right)\right]\right)^{1/q}
\,.
\end{eqnarray}
Thus,
\begin{eqnarray}
E[E_{\omega,\omega}I_N]&\leq&\sum_{t=1}^N E[G(t)]\leq
CN^{2/r'}N^{1/p}\left(\sum_{t=1}^N
E\left[P_{\omega\times\omega}\left(G(t)\neq 0\right)\right]\right)^{1/q}
\end{eqnarray}
%\sqrt{E[G(t)^2]E[{\bf 1}_{G(t)\neq 0}^2]}\\
% &\leq& CN^{2/r} \cdot\sqrt{E\left[P_{\omega\times\omega}\left(G(t)\neq
% 0\right)\right]}\,.
%\end{eqnarray*}
Using (\ref{eq:ren2}), we see that
\[
E[E_{\omega,\omega}I_N]\leq
C N^{\frac{2}{r'}+\frac1p+(\frac{p-1}{p})(\frac12+\frac1r'+a_2)}\,.
%\frac 34 + \frac {a_2}{2} + \frac{5}{2r}}.
\]
By choosing $2p<r'$ close to $r'$ and $a_2$ small, we can get in the last
exponent any power strictly larger than $4/r'+1/2-2/(r')^2$.

\end{proof}
\section*{Acknowledgment}
We thank N.~Zygouras for useful discussions. Section
\ref{addendum}
was written following a very useful conversation with
F. Rassoul-Agha and
T. Sepp\"{a}l\"{a}inen, who  described to one of us
their work in \cite{RS07}. We thank
T. Sepp\"{a}l\"{a}inen for very useful
comments on an earlier draft of this paper.

\end{document}